\documentclass[10pt,twoside,english,reqno,a4paper]{amsart}
\usepackage{listings,graphicx,amsmath,varioref,amscd,amssymb,color,bm,stmaryrd,amsthm,amsfonts,graphics,geometry,latexsym,pgf,pst-all}

\usepackage{color}

\theoremstyle{plain}
\usepackage{esint}
\usepackage{cancel, ulem}
\usepackage{amsthm}

\usepackage{graphicx,fancybox,color,epic,eepic}
\usepackage{ulem}
\usepackage{tikz-cd}
\usepackage{geometry}
\usepackage{graphicx}

\definecolor{red}{rgb}{0.9,0,0}

\definecolor{azul}{rgb}{0,0,1}

\usepackage{hyperref}
\usepackage{tikz-cd}
\usetikzlibrary{patterns}

\usepackage{indentfirst}
\usepackage{graphicx}

\theoremstyle{plain}
\newtheorem{theorem}{Theorem}[section]
\newtheorem{proposition}[theorem]{Proposition}

\newtheorem{lemma}[theorem]{Lemma}

\newtheorem{corollary}[theorem]{Corollary}

\theoremstyle{definition}
\newtheorem{definition}[theorem]{Definition}
\newtheorem{remark}[theorem]{Remark}

\theoremstyle{remark}

\numberwithin{equation}{section}
\def\fd{\hfill{\vbox to 7pt{\hbox to 7pt{\vrule height 7pt width 7pt}}}}

\def\:{\colon}

\def\R{\mathbb R}
\def\z{\mathbf z}
\DeclareMathOperator{\sgn}{sgn}

\usepackage{palatino}
\def\rn{\mathbb{R}^N}

\makeatletter
\@namedef{subjclassname@2020}{  \textup{2020} Mathematics Subject Classification}
\makeatother
\usepackage{geometry}
\title{Quasilinear elliptic equations involving asymptotically  linear operators}

\author[F. Oliva]{Francescantonio Oliva}
\author[F. Petitta]{Francesco Petitta}
\author[M. T. O. Pimenta]{Marcos Tadeu de Oliveira Pimenta}
\author[S. Segura de Le\'on]{Sergio Segura de Le\'on}

\address{Francescantonio Oliva
	\hfill \break\indent
Dipartimento di Scienze di Base e Applicate per l' Ingegneria, Sapienza Universit\`a di Roma
	\hfill \break\indent
	Via Scarpa 16, 00161 Roma, Italy}
\email{\tt francescantonio.oliva@uniroma1.it}
\address{Francesco Petitta 	\hfill \break\indent Dipartimento di Scienze di Base e Applicate per l' Ingegneria, Sapienza Universit\`a di Roma
		\hfill \break\indent
		Via Scarpa 16, 00161 Roma, Italy}
\email{\tt francesco.petitta@uniroma1.it}
\address{Marcos Tadeu de Oliveira Pimenta
	\hfill \break\indent Departamento de Matem\'atica e Computa\c{c}\~ao, Universidade Estadual Paulista - Unesp \hfill\break\indent 19060-900, Presidente Prudente - SP, Brazil.} \email{{\tt
		marcos.pimenta@unesp.br}}
\address{Sergio Segura de Le\'on
	\hfill \break\indent Departament d'An\`alisi Matem\`atica,
	Universitat de Val\`encia, \hfill\break\indent Dr. Moliner 50,
	46100 Burjassot, Val\`encia, Spain.} \email{{\tt
		sergio.segura@uv.es}}

\date{}

\begin{document}

\maketitle
\begin{abstract}
In this paper 
we aim to study in a unified fashion  existence and regularity properties of solutions to a very general class of Dirichlet boundary value problems such as
\begin{equation}\label{P1abs}\tag{1}
	\begin{cases}
		- \mbox{div}\left(\displaystyle \phi(|\nabla u|) \nabla u\right)  =  f(x) &\quad \mbox{in } \Omega,\\
		u =  0 &\quad \mbox{on } \partial \Omega,
	\end{cases}
\end{equation} 
where $\phi:[0,\infty)\to [0,\infty)$ is  any function such that  $\phi(s)s \in C((0,\infty))$ is a non-decreasing function which tends to $1$ as $s\to\infty$. Among  others, most of the  far-famed operators derived from  functionals with linear growth, as for instance the $1$-Laplacian and  the minimal surface operator are included in our model. Under suitable assumptions on the data we prove existence and regularity results for problem \eqref{P1abs} which recover the sharp thresholds known in the model cases.
\end{abstract}
\tableofcontents


\section{Introduction}

Since the 1970s, significant efforts have been devoted by many mathematicians to placing quasilinear elliptic problems involving non-homogeneous divergence form operators of the type 
\begin{equation}\label{P0intro}
	\begin{cases}
		- \mbox{div}\left(\displaystyle \phi(|\nabla u|) \nabla u\right)  =  f(x) &\quad \mbox{in } \Omega,\\
		u =  0 &\quad \mbox{on } \partial \Omega,
	\end{cases}
\end{equation}
 within a unified and coherent theoretical framework
 thus extending the case of power type growth, i.e., $\phi(s)\sim s^{p-2}$.

As far as we know, such studies started with Donaldson in \cite{Donaldson} and Gossez in \cite{Gossez}, giving rise to plenty of works involving what has started to be called the $\Phi-$Laplacian operator, defined as
$$
\Delta_\Phi u := \mbox{div}\left(\displaystyle \phi(|\nabla u|) \nabla u\right),
$$
being $\Phi$ the primitive of $\phi(s)s$ such that $\Phi(0)=0$.

 The  case of $\Phi-$Laplacian for problems as in \eqref{P0intro} has been largely developed since then; without the aim to be exhaustive let us only mention the gentle survey \cite{chl} and references therein. 

Other  relevant studies  concerning $\Phi-$Laplace type operators, mainly related to the case of a reaction type nonlinearity, include the one of   Cl\'ement, Garc\'ia-Huidobro, Man\'asevich and Schmitt \cite{ClementHuidobroManasevichSchmitt}, Garc\'ia-Huidobro, Le, Man\'asevich and Schmitt \cite{HuidobroLeManasevichSchmitt}, Gossez and Man\`asevich \cite{GossezManasevich}, Mih\u{a}ilescu and Repov\v{s} \cite{MihailescuRepovs}. Also, some recent works on this issue are worth mentioning, as \cite{CarvalhoGoncalvesSilva}, where Carvalho, Gon\c{c}alves and Silva established a kind of Br\'ezis-Oswald problem to such type of operator. In \cite{AlvesSantosZhou}, Alves, Abrantes Santos and Zhou studied a semipositone problem involving the $\Phi-$Laplacian operator, while in \cite{SantosSoares}, Abrantes Santos and Soares studied an optimization problem with volume constraint involving the $\Phi-$Laplacian operator.

In all of the papers cited above, the function $s \mapsto s\phi(s)$ is an increasing homeomorphism from $[0,\infty)$ into itself, such that it defines a convex function $\Phi(s) = \int_0^{|s|} \sigma\phi(\sigma) d\sigma$, which is a {\it Young function}, i.e., it satisfies 
$$
\Phi(0) = 0 \quad \mbox{and} \quad \lim_{|s| \to \infty} \Phi(s) = \infty.
$$
Actually, it satisfies stronger conditions:
\begin{equation}
\label{eq:Nfunction}
\lim_{s \to 0} \frac{\Phi(s)}{s} = 0 \quad \mbox{and} \quad \lim_{|s| \to \infty} \frac{\Phi(s)}{s} = \infty,
\end{equation}
which allow us to say that $\Phi$ is an $N-$function (following \cite{Adams} for instance). 

Problems involving the $\Phi-$Laplacian operator are interesting, since they present a unified way to deal with many quasilinear elliptic problems. For instance, 
\begin{itemize}
\item if $\Phi(s) = \frac{s^p}{p}$ and $1 < p $, we have the $p$-Laplacian operator;
\item if $\Phi(s) = \frac{s^p}{p} + \frac{s^q}{q}$, where $1 < p < q $, we have the $(p,q)-$Laplacian operator;
\item if $\Phi(s) = s^p \log (1 + s)$ with $p\ge 1$, we have a plasticity type problem.
\end{itemize}
A crucial starting point in the approach followed by all these papers, is the fact that \eqref{eq:Nfunction} allows them to set the problem \eqref{P0intro} in the Orlicz-Sobolev space $W^{1,\Phi}_0(\Omega)$. Its definition is a bit technical (we refer to \cite{Adams}), but roughly speaking, it is the space of functions in the Orlicz space $L_\Phi(\Omega)$, such that all their first order derivatives also belong to this space. The Orlicz space $L_\Phi(\Omega)$, in turn, is formed by those measurable functions $u: \Omega \to \mathbb{R}$, such that
$$
\|u\|_\Phi := \sup \left\{ \int_\Omega u v dx: \, \int_\Omega \Phi^*(|v|) dx \leq 1\right\} < \infty,
$$
where
$$
\Phi^*(t) = \sup \{st - \Phi(s); \, s \geq 0\}
$$
is the {\it convex conjugate} of $\Phi$ (or the {\it complementary function} of $\Phi$).

It is worthwhile to observe that several important quasilinear elliptic equations can be expressed as in \eqref{P0intro}, but with associated functions $\phi(s)s$ which give rise to functions $\Phi$ that do not satisfy \eqref{eq:Nfunction}. Hence, in these situations, a treatment based on setting the problem in Orlicz-Sobolev spaces as handled in \cite{Adams} is no longer possible. As examples, we can cite the problems involving the mean-curvature operator as
\begin{equation}
\label{P1MC}
	\begin{cases}
		- \mbox{div}\, \left(\displaystyle \frac{D u}{\sqrt{1 + |D u|^2}}\right)  =  f(x) &\quad \mbox{in } \Omega,\\
		u =  0 &\quad \mbox{on } \partial \Omega,
	\end{cases}
\end{equation}
and also those ones involving the $1$-Laplacian operator
\begin{equation}
\label{P1laplacian}
	\begin{cases}
		- \mbox{div}\, \left(\displaystyle \frac{D u}{|D u|}\right)  =  f(x) &\quad \mbox{in } \Omega,\\
		u =  0 &\quad \mbox{on } \partial \Omega.
	\end{cases}
\end{equation}
Indeed, in the former, $\Phi(s) = \sqrt{1 + s^2} - 1$, while in the latter, $\Phi(s) = |s|$.  It is worth underlining that, in contrast with \eqref{eq:Nfunction}, for both problems it holds $\lim_{s \to \infty} \displaystyle \frac{\Phi(s)}{s} = 1$ and moreover, for the latter $\lim_{s \to 0} \displaystyle \frac{\Phi(s)}{s}$ does not exist. In both cases, it is known that if the datum $f$ belongs to $L^N(\Omega)$, then a solution exists only when $\| f\|_{L^{N}(\Omega)}$ is small enough (see \cite{gia,g76,g78, gop} and  \cite{K1,K2,CT, MST}  respectively, and references therein).

Taking all of these questions into account, 
in this paper, we aim to study the general problem
\begin{equation}\label{P1intro}
	\begin{cases}
		- \mbox{div}\, \left(\displaystyle \phi(|D u|) D u\right)  =  f(x) &\quad \mbox{in } \Omega,\\
		u =  0 &\quad \mbox{on } \partial \Omega,
	\end{cases}
\end{equation}
where $\Omega \subset \mathbb{R}^N$ is a bounded open set having Lipschitz continuous boundary, $N \geq 2$ and $f \in L^N(\Omega)$. The function $\phi:[0,\infty)\to [0,\infty)$ is assumed to satisfy: 
\begin{equation}\label{phi1}
	\phi(s) s \in C((0,\infty)) \text{ and }   \phi(s)s=0 \text{ at } s=0;
\end{equation}
\begin{equation}\label{phi2}
	\phi(s) s \text{ is non-decreasing for}\ s\in [0,\infty);
\end{equation}
\begin{equation}\label{phi3}
	\lim_{s \to \infty}\phi(s)s = 1.
\end{equation}
Notice that \eqref{phi2}  implies that there exists $ \lim_{s\to 0^+} \phi(s)s$ and it might happen that this limit is strictly positive.

The approach we follow in this work, is strongly influenced by the approximation scheme established in the 2000's for problems involving the $1$-Laplacian operator (see \cite{AndreuBallesterCasellesMazon,AndreuBallesterCasellesMazonJFA,D}), where the original equation is approximated by $p$-Laplacian problems, uniform bounds on the solutions are obtained and then the limit is taken in some appropriate sense. 
A notion of solution to problem \eqref{P1laplacian} is introduced in \cite{AndreuBallesterCasellesMazon,D} as follows:
 $u \in BV(\Omega)$ is a solution if there exists ${\textbf z} \in L^\infty(\Omega,\mathbb{R}^N)$ such that $\mbox{div}\, {\textbf z} \in L^N(\Omega)$, $\|{\textbf z}\|_\infty \leq 1$ and
\begin{equation*}\label{weak1}
\begin{cases}
- \mbox{div}\, \z = f(x) & \mbox{ in } \mathcal{D}'(\Omega),\\
(\z, Du)  = |Du| & \mbox{ in } \mathcal{M}(\Omega),\\
[\z,\nu] \in \sgn(-u) &\ \mathcal{H}^{N-1}\text{-a.e. on } \partial \Omega.
\end{cases}
\end{equation*}
In the last expression, the second equation is intrinsically related to the $1$-Laplacian operator and implies that in a generalized sense, ${\bf z}$ (along with condition $\| \z\|_\infty\le 1$) plays the role of $D u/|D u|$ in \eqref{P1laplacian}. 

One basic question which motivates our work is:
{\it What conditions should be imposed on $\z$ and $Du$, in order to guarantee that $\z$ plays the role of $\phi(|Du|)Du$?}

An answer to this question, for the special case $\phi(s)s = s/\sqrt{1+s^2}$, has been given in \cite{Te,dt} (see also \cite{stevsch,SchevenSchmidtAIHP,gop}). More specifically, their notion of solution contains the condition that the vector field $\mathbf{z}$, playing the role of $\displaystyle \frac{Du}{\sqrt{1+|Du|^2}}$, must satisfy:
\[(\z, Du) =  \sqrt{1 + |Du|^2} - \sqrt{1 - |{\textbf z}|^2}.\]
For more recent related works, we refer to \cite{MontenegroPimenta,OlivaPetittaSegura}.

In this paper, we propose a  unified  approach to deal with \eqref{P1intro}, by giving, as a by-product,  a detailed answer to the question raised above. Our results contain sharply the known results in the case of most of the famous operators including the 1-Laplace operator and  the minimal surface operator.  
The arguments rely strongly on the Fenchel--Young inequality (see Subsection \ref{FYI} and condition \eqref{eqweakform2} in Definition \ref{def} below).

Our main result states that if $f\in L^N(\Omega)$ and $\|f\|_{L^{N}(\Omega)}<\mathcal{S}_{1}^{-1}$, where $\mathcal{S}_1$ denotes the best constant of Sobolev's embedding $W_0^{1,1}(\Omega)\hookrightarrow  L^{\frac N{N-1}}(\Omega)$, then there exists a bounded solution to problem \eqref{P1intro}. The approach we follow relies on the idea of approximating the operator $\displaystyle \mbox{div}\, \left(\displaystyle \phi(|D u|) D u\right)$ by operators having $(p-1)$--growth, for which it can be applied the classical theory of monotone operators in $W_0^{1,p}(\Omega)$, with $p>1$. We finally achieve our goal, by obtaining a priori estimates and letting $p$ go to $1$. 

As far as we know, the only paper which has addressed a related question (but using a different approach), is \cite{AndreuCasellesMazon}, in which, in case ${\phi(s)s} \to  0$ as ${s \to 0^{+}} $, the authors considered an elliptic equation as in \eqref{P1intro} with an absorption term. The existence result in this case is then used in order to get existence of a solution  for a parabolic quasilinear equation  involving general functionals with asymptotically linear growth.

In Section 2 we present some preliminaries concepts and definitions. In Section 3 we perform a deep study of the function $\Phi$, by introducing some concepts of convex analysis and presenting the condition that ${\bf z}$ should satisfy, ensuring that it plays the role of the vector field appearing into the divergence in \eqref{P0intro}. In Section 4 we finally present our main result. In Section 5 we present a condition on $\phi$ under which one can prove that the vector field $\z$ is uniquely determined. In Section 6, we extend the  result  of our main theorem in case of data $f$ in the Lorentz space $L^{N,\infty}(\Omega)$.

\section{Preliminaries}
\label{Sectionpreliminaries}

\subsection{Notation}

For the entire paper, $\Omega$  denotes an open bounded subset of $\R^N$ ($N\ge 2$) with Lipschitz boundary.
By $\mathcal H^{N-1}(E)$ we mean the $(N - 1)$-dimensional Hausdorff measure of a set $E$  while $|E|$ stands for its $N$-dimensional  Lebesgue measure. 

Unless otherwise specified, by $C$ we denote several positive constants whose value may change from line to line and, sometimes, on the same line. These values will only depend on the data but they will never depend on the indexes of the sequences we will gradually introduce. Let us explicitly mention that we will not relabel an extracted compact subsequence.

We denote by ${\rm sgn} (s)$ the multi-valued sign function defined by 
$$
\sgn (s):=\begin{cases}
	1 & \text{if }\ s>0,\\
	[-1,1] & \text{if }\ s=0,\\
	-1 & \text{if }\ s<0.
\end{cases}
$$

In what follows, we use two auxiliary real functions. Given $k>0$, we define $T_k$ and $G_k$ as
\[T_k(s)=\left\{\begin{array}{ll}
s&\hbox{if }|s|\le k,\\
k\, \sgn(s) &\hbox{if } |s|> k,\end{array}\right.\qquad G_k(s)=(|s|-k)^+\sgn(s). \]

As usual, we write $L^q(\Omega)$ and $W_0^{1,p}(\Omega)$ to indicate, respectively, Lebesgue spaces and Sobolev spaces of functions whose trace on the boundary vanishes. 
We denote by $\mathcal{S}_1$ the best constant in the Sobolev inequality for functions in $W^{1,1}_0(\Omega)$, that  is
\begin{equation}\label{sob}
	\|v\|_{L^{1^*}(\Omega)} \le \mathcal{S}_1 \|v\|_{W^{1,1}_0(\Omega)}, \ \ \forall v \in W^{1,1}_0(\Omega),\end{equation}
where $1^*=\frac{N}{N-1}$. 

\subsection{Functions of bounded variation}

The natural space to deal with \eqref{P1intro} is the space of functions of bounded variation defined as
$$BV(\Omega):=\{ u\in L^1(\Omega) : Du \in \mathcal{M}(\Omega,\rn)\}.$$  
Here $Du$ denotes the distributional gradient of $u$ and by
 $Du \in \mathcal{M}(\Omega,\rn)$ we mean that  each distributional partial derivative of $u$ is a Radon measure  having finite total variation over $\Omega$. The total variation of $Du$ over $\Omega$ is denoted by $\int_\Omega |Du|$. 

We remark that every BV-function exhibits a trace on $\partial\Omega$, so that $u\big|_{\partial\Omega}$ has a meaning.
We underline that  the $BV(\Omega)$ space endowed with the norm  
$$\displaystyle \|u\|_{BV(\Omega)}=\int_{\partial\Omega}
|u|\, d\mathcal H^{N-1}+ \int_\Omega|Du|,$$
is a Banach space. Just like $W_0^{1,1}(\Omega)$, this Banach space is continuously embedded in $L^{1^*}(\Omega)$ and the embedding of $BV(\Omega)$ in each $L^q(\Omega)$ with $1\le q<1^*$ is compact.

We refer to \cite{AmbrosioFuscoPallara} for a complete account on $BV$-spaces and for further standard notations not mentioned here.

\subsection{The Anzellotti-Chen-Frid theory}\label{sec:anz} Here we briefly present the $L^\infty$-divergence-measure vector fields theory due to \cite{An} and \cite{CF}. We denote by 
$$
X_N(\Omega) = \{\z \in L^\infty(\Omega,\mathbb{R}^N); \, \operatorname{div}\z \in L^N(\Omega)\}.
$$
In \cite{An}, the following distribution $(\z,Dv): C^1_c(\Omega)\to \mathbb{R}$ is defined: 
\begin{equation*}\label{dist1}
	\langle(\z,Dv),\varphi\rangle:=-\int_\Omega v\varphi\operatorname{div}\z-\int_\Omega
	v\z\cdot\nabla\varphi,\quad \varphi\in C_c^1(\Omega).
\end{equation*}
Observe that if $v\in BV(\Omega)$ and $\z\in X_N(\Omega)$ then the above distribution is well defined.
Moreover, it is a Radon measure and it holds 
\begin{equation*}\label{finitetotal1}
	\left| \int_B (\z, Dv) \right|  \le  \int_B \left|(\z, Dv)\right| \le  \|\z\|_{L^\infty(U,\rn)} \int_{B} |Dv|\,,
\end{equation*}
for all Borel sets $B$ and for all open sets $U$ such that $B\subset U \subset \Omega$.
We recall that  every $\z \in X_N(\Omega)$ possesses a weak trace on $\partial \Omega$ of its normal component which is denoted by
$[\z, \nu]$, where $\nu(x)$ is the outward normal unit vector defined for $\mathcal H^{N-1}$-almost every $x\in\partial\Omega$ (see  \cite{An}). Moreover, it holds
\begin{equation*}\label{des1}
	\|[\z,\nu]\|_{L^\infty(\partial\Omega)}\le ||\z||_{L^\infty(\Omega,\rn)}\,.
\end{equation*}

Furthermore, if $\z\in X_N(\Omega)$ and $v\in BV(\Omega)$, then the following Green formula holds:
	\begin{equation}\label{green}
		\int_{\Omega} v \operatorname{div}\z \, dx+ \int_{\Omega} (\z, Dv) = \int_{\partial \Omega} v[\z, \nu] \ d\mathcal H^{N-1}.
	\end{equation}

\section{Some basic remarks on the  function $\phi$}\label{analysis}
Along the entire paper, the function $\phi:[0,\infty)\to [0,\infty)$, which defines the principal operator into \eqref{P1intro}, satisfies assumptions \eqref{phi1}, \eqref{phi2} and \eqref{phi3} above.
In this section we discuss these structural hypotheses as well as relevant  examples of  equations whose principal operator satisfies our assumptions. As  we are interested in proving existence  results  of solutions to \eqref{P1intro}, we will employ the Fenchel--Young inequality which will be useful in order to
 \begin{enumerate}
 \item give a meaning  to $\phi(|Du|)Du$, when $Du$ is just a Radon measure;
 \item consider functions of the gradient $Du$ that are weak lower semicontinuous (a key ingredient in our arguments).
 \end{enumerate}

\subsection{A glimpse on 1-D convex analysis}

Let $\phi$ be a function satisfying \eqref{phi1}, \eqref{phi2} and \eqref{phi3}.
Since the function $s\mapsto \phi(s)s$ is non-decreasing,  a generalized inverse always exists, and it is given by 
\begin{equation}\label{inverse}
\psi(t)=\inf\{s\ge 0\>:\> \phi(s)s\ge t\}\,,\quad t\ge 0\,.\end{equation}
This inverse is non-decreasing, left--continuous and has an asymptote at $t=1$. As a consequence of the existence of this asymptote, $\psi(t)=\infty$ for all $t>1$.

We point out that we assume that $\phi(s)s$ must be continuous on $(0,\infty)$, but without any further smoothness requirements. We also remark that it may be constant on an interval (so that its generalized inverse may be discontinuous, see Figure 1).  Even a kind of devil's staircase is allowed. 
     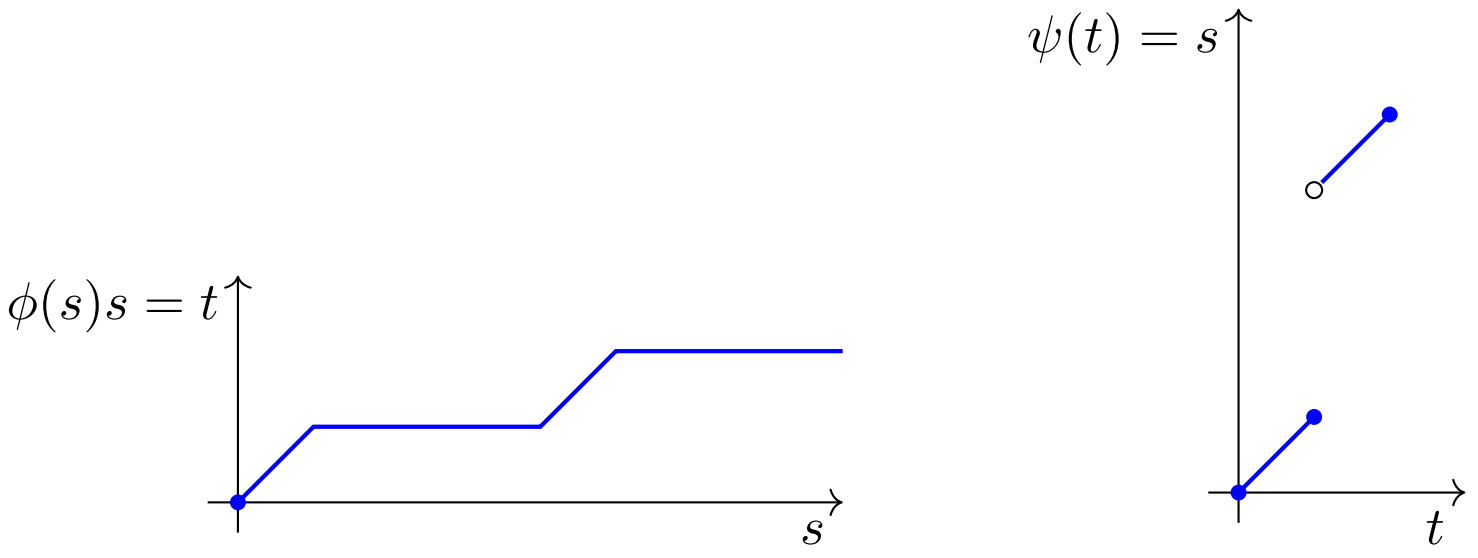
\begin{figure}[h]
\centering
\begin{tikzpicture}
\draw [->,thin] (-0.2,0) -- (4,0);
\draw [->,thin] (0,-0.2) -- (0,1.5);
\draw [thick, blue] (0,0)--(0.5,0.5)--(2,0.5)--(2.5,1)--(4,1);
\fill [blue] (0,0) circle (1.5pt);
\node [left] at (0,1.3) { $\phi(s)s=t$};
\node [below] at (3.8,0) { $s$};
\end{tikzpicture}
\hskip1cm
\centering
\begin{tikzpicture}
\draw [->,thin] (-0.2,0) -- (1.5,0);
\draw [->,thin] (0,-0.2) -- (0,3.2);
\draw [thick, blue] (0,0)--(0.5,0.5);
\draw [thick, blue] (0.55,2.05)--(1,2.5);
\fill [blue] (0,0) circle (1.5pt);
\fill [blue] (0.5,0.5) circle (1.5pt);
\fill [blue] (1,2.5) circle (1.5pt);
\draw (0.5,2) circle (1.5pt);
\node [below] at (1.3,0) { $t$};
\node [left] at (0,3) { $\psi(t)=s$};
\end{tikzpicture}
\caption{Example of $\phi(s)s$ and its inverse}
%\label{,,}
\end{figure}

We consider the function  $\Phi: \mathbb{R} \to \mathbb{R}$ given by 
\begin{equation}
	\label{Phi1}
	\Phi (s) = \int_0^{|s|} \phi (\sigma) \sigma\, d\sigma,
\end{equation}
which is convex and $C^1(\mathbb{R}\backslash\{0\})$. Moreover, we define 
\begin{equation*}
	\label{Psi1}
	\Phi^* (t) = \int_0^{|t|} \psi (\sigma) \, d\sigma,
\end{equation*}
which is its convex conjugate. As pointed out in \cite[Chapter 8]{Adams}, this function could be alternatively defined via Legendre transformation, i.e. as 
$$
\Phi^*(t) = \sup \{st - \Phi(s); \, s \geq 0\}.
$$
Observe that $\Phi^* (t) =\infty$ for every $|t|>1$. It is well--known that the inequality 
\begin{equation*}\label{fenchel_dim_1}
	st\le \Phi(s)+\Phi^*(t)
\end{equation*}
holds for $s,t\ge0$, with equality if and only if $t \in \partial \Phi(s)$, where $\partial \Phi(s)$ denotes the subdifferential of $\Phi$ in $s \in \mathbb{R}$. 

Actually, it is a  generalization of the usual Young inequality, which follows from it by consi\-de\-ring $\Phi(s) = \frac{1}{q}|s|^q$, where $q > 1$; here $\Phi^*(s) = \frac{1}{q'}|s|^{q'}$, with $\frac{1}{q}+ \frac{1}{q'}=1$. 

 For details and properties of functions as $\Phi$ (and $\Phi^*$) we refer to \cite{BS} and \cite{Z}.

\subsection{Fenchel--Young's inequality}\label{FYI}

Here we recall some results from convex analysis in $\mathbb{R}^N$, $N>1$ (see \cite{EkelandTeman} for a complete reference). Let   $\Gamma:\mathbb{R}^{N}\to\overline{\mathbb{R}}$ be a convex function. Then the function $\Gamma^{*}:\mathbb{R}^{N}\to\overline{\mathbb{R}}$ given by
\begin{equation*}
\label{Convexconjugate}
\Gamma^{*}(\eta)=\sup\{\xi\cdot \eta - \Gamma(\xi): \xi\in\mathbb{R}^{N}\}
\end{equation*}
is called the convex conjugate of $\Gamma$ and it can be proven to be convex. From the definition of $\Gamma^*$, it holds the so called Fenchel--Young inequality:
\begin{equation}
\label{YFinequality}
\Gamma(\xi) + \Gamma^{*}(\eta)\geq \xi\cdot \eta \qquad \text{for all }\xi, \eta\in\mathbb{R}^{N},
\end{equation}
with equality if and only if $\eta \in \partial \Gamma(\xi)$.

Hence, for every $\xi \in \mathbb{R}^N$ where $\Gamma$ is $C^1$, the equality in \eqref{YFinequality} will be satisfied at $\eta \in \mathbb{R}^N$, such that
\begin{equation}
\label{YFequality}
\eta = \nabla \Gamma(\xi).
\end{equation}
 
Then
\begin{equation}
\label{eqConvexconjugate}
\Gamma^*(\eta) = \xi \cdot \nabla \Gamma(\xi) - \Gamma(\xi),
\end{equation}
where $\eta = \nabla \Gamma(\xi)$.

In order to be concrete, let's briefly discuss how we will apply the previous argument.

We want to apply \eqref{eqConvexconjugate} in order to fully understand the principal operator defined in \eqref{P1intro}.
With this goal in mind, the convex function we are interested in is
 $\Gamma: \mathbb{R}^N \to \mathbb{R}$, given by $\Gamma(\xi)=\Phi(|\xi|)$ (with $\Phi$ defined in \eqref{Phi1} above). 
 Then, for $u \in W^{1,1}_0(\Omega)$ and $\z \in L^\infty(\Omega,\mathbb{R}^N)$, such that $\|\z\|_\infty \leq 1$, by \eqref{eqConvexconjugate} and \eqref{YFequality}, it follows that
  \begin{equation} 	\label{YoungFencheleq1BIS}
 	\z\cdot \nabla u \le \Phi(|\nabla u|) + \Phi^*(|\z|),
 \end{equation}
where
 \begin{equation} 	\label{YoungFencheleq1}
 	\z\cdot \nabla u = \Phi(|\nabla u|) + \Phi^*(|\z|) \Leftrightarrow \z = \phi(|\nabla u|)\nabla u,
 \end{equation}
 whenever $\nabla u \neq 0$.

 Therefore the equality on the left-hand side of the previous equivalence  can be seen as a  way to interpret  vector field $\phi(|\nabla u|)\nabla u$ that appears in problem \eqref{P1intro}. To be more precise, as the solutions to \eqref{P1intro} are, in general, only $BV$-functions (i.e. $D u$ is only a vectorial measure) we need a relaxed version of \eqref{YoungFencheleq1} which may involve equalities between measures.  On the left-hand side, the dot product has to be replaced by the Anzellotti pairing $(\z, Du)$, which is well--defined as long as $u\in BV(\Omega)$ and $\z\in L^\infty(\Omega, \R^N)$ is such that 
$\hbox{div}\, \z \in L^N(\Omega)$. 

 We need to define $\Phi(|Du|)$ where  $Du$ is a Radon measure. 
First observe that the recession function of $\Gamma(\xi)$, i.e. $\Gamma^0(\xi) $,  in this case  is defined by
\begin{equation}\label{super0}
\Gamma^0(\xi) = \lim_{t \to \infty} \frac{\Phi(t|\xi|)}{t}=|\xi|.
\end{equation}
Now, let $\mu \in {\mathcal{M}(\Omega,\rn)}$ and let $\mu = \mu^a + \mu^s$ be its Radon-Nikod\'ym decomposition with respect to the $N-$dimensional Lebesgue measure $\mathcal{L}^N$. Having in mind \eqref{super0} (see also \cite[Definition $2.3$]{anzellotti}), it is natural to define the measure $\Phi(|\mu|) = \Gamma (\mu)$ as
$$
\int_B \Phi(|\mu|) := \int_B \Phi(|\mu^a|)\, dx + \int_B|\mu^s|\,,
$$
for all Borel sets $B\subset\Omega$. In case $u\in BV(\Omega)$, we obtain
\begin{equation}\label{defmeasure}
\int_B \Phi(|Du|) := \int_B \Phi(|D^au|)\, dx + \int_B|D^su|\,,
\end{equation}
for all Borel sets $B\subset\Omega$.
  Therefore, a relaxed version of \eqref{YFinequality} still holds, namely: 
  let $u \in BV(\Omega)$ and $\z \in L^\infty(\Omega,\mathbb{R}^N)$ such that 
   $\hbox{div}\,\z\in L^N(\Omega)$ and $\|\z\|_\infty \leq 1$, then
 \begin{equation}\label{YoungFencheleq2}
 (\z,D u)\le \Phi(|D u|) + \Phi^*(|\z|)\mathcal{L}^N \quad \mbox{in $\mathcal{M}(\Omega)$}.
 \end{equation}  
 Indeed, splitting the measure $(\z,Du)$ into its absolutely continuous part and its singular part with respect to the Lebesgue measure and recalling also \cite[Theorem $2.4$]{An}, it yields 
 $$
\begin{aligned}
 	(\z, Du) & =  (\z, Du)^a + (\z, Du)^s
 	=  \z\cdot D^au + (\z, D^su)\\
 	& \leq  \Phi(|Du^a|)\mathcal{L}^N + |D^su| + \Phi^*(|\z|)\mathcal{L}^N\\
 	& =  \Phi(|Du|) + \Phi^*(|\z|)\mathcal{L}^N,
 \end{aligned}
 $$
 where we made advantage of \eqref{YoungFencheleq1BIS} and of \eqref{defmeasure}.
 Moreover, in accordance with \eqref{YoungFencheleq1}, we mean that $\z = \phi(|D u|)D u$ once it holds
 \begin{equation}\label{YoungFenchelIdweak}
 	(\z,D u) = \Phi(|D u|) + \Phi^*(|\z|)\mathcal{L}^N,
 \end{equation}  
 which is equivalent to have 
 \begin{equation}\label{YoungFenchelIdweak1}
 	 (\z,Du)^a =  \z\cdot D^a u= \Phi(|D^a u|) + \Phi^*(|\z|)
 \end{equation}
 and 
 \begin{equation}\label{YoungFenchelIdweak2}
 	 (\z,Du)^s = (\z,D^s u) = \Phi(|D^s u|) =|D^s u|,
 \end{equation}
 as, once again, the first equality in \eqref{YoungFenchelIdweak1} follows from \cite[Theorem $2.4$]{An}. 

\subsection{Lower semicontinuity results}
Now let us state some technical results, dealing with the lower semicontinuity of suitable functionals.

\begin{lemma}
	\label{lemma.lsc1}
	For $0\leq \varphi \in \mathcal{D}(\Omega)$, let $w_n ,w\in BV(\Omega)$ be such that $w_n \to w$ in $L^1(\Omega)$. Then
	\begin{equation}
		\label{eq.lsc1}
		\int_\Omega \varphi \Phi(|Dw|) \leq \liminf_{n \to \infty} \int_\Omega \varphi \Phi (|D w_n|)
	\end{equation}
	and
	\begin{equation}
		\label{eq.lsc2}
		\int_\Omega \Phi(|Dw|) + \int_{\partial \Omega}|w| d\mathcal{H}^{N-1} \leq \liminf_{n \to \infty} \int_\Omega \Phi(|D w_n|)+ \int_{\partial \Omega}|w_n| d\mathcal{H}^{N-1},
	\end{equation}
 where $\Phi$ is defined in \eqref{Phi1}.	
\end{lemma}
\begin{proof}
%Let us prove \eqref{eq.lsc2}, as \eqref{eq.lsc1} is similar. 
As $\Phi$ is convex and continuous, the proof follows from \cite[Fact 3.4]{anzellotti}.
\end{proof}

\subsection{Some relevant examples of operators}\label{sec_genop}

We start supplying examples of functions $\phi (s)$, the corresponding functions  $\Phi(s)$ and $\Phi^*(t)$, and sometimes the associated operator ${\rm div}  ( \phi(|Du|)Du)= {\rm div}  (\nabla_\xi \Phi(|Du|))$:
	
	\begin{enumerate}
	\item If $\phi(s)s=1$ for all $s>0$, then the principal operator is given by the $1$-Laplacian. In this case, equation \eqref{P1intro} reads as $-\Delta_1 u=f(x)$. The associated convex functions are given by $\Phi(s)=|s|$ and $\Phi^*(t)=0$ if $|t|\le 1$.
		\item If $\phi(s)s = \frac{s}{\sqrt{1+s^2}}$, then in equation \eqref{P1intro} one gets the minimal surface operator. Here we have $\Phi(s)=\sqrt{1+s^2}-1$ and $\Phi^{*}(t)=1-\sqrt{1-t^2}$ when $|t|\le 1$.
		\item Let 
		\[
		\phi(s)s=
		\begin{cases}
			s & \hbox{ if }0\le s\le 1,\\
			1 & \hbox{ if }s>1.
		\end{cases}
		\]
		Then we obtain
		\[
		\Phi(s)=
		\begin{cases} \displaystyle 
			\frac12 s^2 &\hbox{ if }  |s|\le 1, \\
			|s|-\frac12 &\hbox{ if } |s|>1.
		\end{cases}
		\]
		and $\Phi^{*}(t)=\frac12 t^2$ when $|t|\le 1$.
		In this case, the associated operator arises when studying plastic antiplanar shear (see \cite{Zh}).
	\end{enumerate}
	
	Other simple examples are given by
	
	\begin{enumerate}
		\item Let $\phi(s) s=\frac s{1+s}$. It follows that $\Phi(s)=|s|-\log(1+|s|)$ and $\Phi^{*}(t)=-|t|-\log(1-|t|)$  when $|t|\le 1$.
		\item If $\phi(s)s=\frac 2\pi \arctan s$, then $\Phi(s)=\frac1\pi(2s\arctan s-\log(1+s^2))$ and $\Phi^*(t)=-\frac2\pi\log(\cos(\frac\pi2 t))$ for $|t|\le 1$. 
		\item Let $\phi(s)s=T_1\left(G_{\frac{1}{2}}(s)\right)$ for $s\ge 0$. Then 
		\[
		\Phi(s)= 
		\begin{cases}\displaystyle
			0 &\hbox{ if } |s|\le \frac12,\\
			\displaystyle \frac{s^2}{2}-\frac{|s|}{2}+\frac18 &\hbox{ if } \frac12< |s|\le  \frac{3}{2},\\
			\displaystyle|s|- 1 &\hbox{ if } |s|> \frac{3}{2},
		\end{cases}
		\] 
		and $\Phi^*(t) = \frac{t^2}{2} + \frac{|t|}{2}$ once $|t|\le 1$. 
		
	\end{enumerate}

\section{Main assumptions, existence and regularity results}

Consider problem 
\begin{equation}
 \label{P1}
\begin{cases}
- \mbox{div}\, \left(\displaystyle \phi(|D u|) D u\right) = f(x) \quad &\mbox{in } \Omega,\\
u  =  0 \quad &\mbox{on } \partial \Omega,
 \end{cases}
\end{equation}
where $\phi$ satisfies \eqref{phi1}, \eqref{phi2} and \eqref{phi3} and $f \in L^N(\Omega)$.

Motivated by the discussion in Section  \ref{FYI} we set what we mean by a solution to \eqref{P1}.

\begin{definition}\label{def}
A function  $u \in BV(\Omega)$ is a  solution of \eqref{P1}, if there exists $\z \in X_N(\Omega)$, such that $\|\z\|_\infty \leq 1$ and
\begin{equation}
-\operatorname{div}\z = f(x) \quad \mbox{in $\mathcal{D}'(\Omega)$},\label{eqweakform1}
\end{equation}
\begin{equation}
(\z,Du) = \Phi(|Du|) + \Phi^{*}(|\z|)\mathcal{L}^N  \quad \mbox{in $\mathcal{M}(\Omega)$},\label{eqweakform2}
\end{equation}
\begin{equation}
{[\z,\nu]}  \in \mbox{${\rm sgn}(-u)$ $\mathcal{H}^{N-1}$-a.e. on $\partial \Omega$}.\label{eqweakform3}
\end{equation}

\end{definition}
\begin{remark}\label{observacion}
\rm
A few words are needed regarding Definition \ref{def}; firstly, it follows from the discussion involving \eqref{YoungFenchelIdweak}, \eqref{YoungFenchelIdweak1} and \eqref{YoungFenchelIdweak2} that \eqref{eqweakform2} is the weak way in which $\z$ plays the role of $\phi(| D u|) D u$ (see Section \ref{s5} for further insight on this identification). Then, \eqref{eqweakform1} and \eqref{eqweakform2} represent how the first equation in \eqref{P1} is intended. The boundary condition is given by \eqref{eqweakform3}, which is nowadays classical for $BV$ solutions to equations deriving from functionals with linear growth (see for instance \cite{AndreuBallesterCasellesMazonJFA} and references therein).
     
It is also worth mentioning that condition \eqref{eqweakform2} is the natural extension of the analogous one in the relevant explicit cases of the mean curvature operator and of the $1$-Laplacian. Specifically, in the case of the mean curvature operator, $\phi(s)s=\frac{s}{\sqrt{1+s^2}}$, which gives that \eqref{eqweakform2} reads as $$(\z,Du)=\sqrt{1+|Du|^2} - \sqrt{1-|\z|^2},   \quad \mbox{in $\mathcal{M}(\Omega)$}$$ 
(see \cite{dt,AndreuCasellesMazon}, see also \cite{gop,SchevenSchmidtAIHP,MontenegroPimenta} and references therein), 
while for the $1$-Laplacian one gets $\phi(s)s=1$ and condition \eqref{eqweakform2} yields to 
\[(\z,Du)=|Du|,   \quad \mbox{in $\mathcal{M}(\Omega)$}\]
(see \cite{AndreuBallesterCasellesMazon, AndreuBallesterCasellesMazonJFA, D}).
\end{remark}

The main result of this section  is stated as follows.

\begin{theorem}\label{theorem1}
Let $\phi$ satisfy \eqref{phi1}, \eqref{phi2} and \eqref{phi3} and let $f \in L^N(\Omega)$ be such that
\begin{equation}\label{hf}
	\|f\|_{L^{N}(\Omega)}<\mathcal{S}_{1}^{-1}.
\end{equation} 
Then there exists  a bounded solution of \eqref{P1}.
\end{theorem}
\begin{remark}\rm
	One may wonder what happens if requirement \eqref{phi3} is replaced with $\lim_{s \to \infty}\phi(s)s = A$ for certain $A>0$. Applying Theorem \ref{theorem1} itself, it is straightforward that the smallness condition to get existence becomes $\|f\|_{L^{N}(\Omega)}<A\mathcal{S}_{1}^{-1}$. Nevertheless, then the concept of solution must be modified by assuming that the vector field satisfies $\|\z\|_\infty\le A$.
\end{remark}
\subsection{Approximation scheme, a priori estimates and existence of limit function}

We prove Theorem \ref{theorem1} through suitable approximation arguments which allow us to set the problem in $W^{1,p}_0(\Omega)$. 

The main idea is to consider a $p$-Laplacian type operator for $p>1$; let us consider
\begin{equation}
	\label{Pp2}
	\begin{cases}
		- \mbox{div}\left(\phi(|\nabla u_p|) |\nabla u_p|^{p-1}\nabla u_p\right)  =  f(x) \quad  &\mbox{in } \Omega, \\
		u_p  =0 \quad   &\mbox{on } \partial \Omega.
	\end{cases}
\end{equation}

We are formally taking advantage of the above approximation as $\phi(s)s^{p}$ is continuous and tends, as $p\to 1^+$, to $\phi(s)s^{}$. 

Let us also stress that, for \eqref{Pp2}, all conditions for applying the theory of monotone operators in $W_0^{1,p}(\Omega)$ are straightforward, except coerciveness, which we show below.

\begin{lemma}\label{coerc}
For every $p>1$, there exist two constants $0<\alpha<1$ and $B_p\ge0$, such that $B_p$ is bounded for $p$ close to $1$, and
\[\alpha|\xi|^p-B_p\le \phi(|\xi|)|\xi|^{p+1}\qquad\forall \xi\in\R^N\,.\]
Moreover, if $\lim_{s\to 0^+}\phi(s)s>0$, then $B_p=0$ for all $p>1$.
\end{lemma}
\begin{proof}
Given $0<\alpha<1$ and having in mind \eqref{phi3}, we may find $h>0$ such that $\phi(s)s>\alpha$ for all $s>h$. 
Since $\alpha s-\phi(s)s^{2}$ defines a continuous function on $[0,h]$, there exists $B\ge 0$ such that 
\[\alpha s-B\le \phi(s)s^{2} \qquad\forall s\in[0, h]\,.\]
Thus, for $p>1$, 
\[\alpha s^p-Bh^{p-1}\le \phi(s)s^{p+1} \qquad\forall s\ge0\,.\]
Hence, it is enough to take $B_p=Bh^{p-1}$.

On the other hand, if $\lim_{s\to 0^+}\phi(s)s>0$, one easily checks that the result is fulfilled with $\alpha=\lim_{s\to 0^+}\phi(s)s$ and $B_p=0$ for all $p>1$. 
\end{proof}

On account of Lemma \ref{coerc}, we may invoke \cite{LL} and \cite{Stampacchia} to find a solution $u_p\in W_0^{1,p}(\Omega)\cap L^\infty(\Omega)$.

Next we show that $u_p$ is uniformly bounded in $BV(\Omega)\cap L^\infty(\Omega)$ with respect to $p$.

\begin{lemma}
\label{lemma.boundedness2}
Let $\phi$ satisfy \eqref{phi1}, \eqref{phi2} and \eqref{phi3} and let $f \in L^N(\Omega)$ be such that \eqref{hf} is in force. Let $u_p$ be a solution to \eqref{Pp2}. Then $u_p$ is uniformly bounded in $BV(\Omega)\cap L^\infty(\Omega)$ with respect to $p$. Moreover it holds
\begin{equation}\label{3222}
\int_\Omega |\nabla u_p|^{p} dx \leq C,
\end{equation}
for some positive constant $C$ independent of $p$.
\end{lemma}
\begin{proof}
 	Firstly we show that $u_p$ is uniformly bounded in $L^\infty(\Omega)$ with respect to $p$. Let $A_k:=\{x\in\Omega: |u_p(x)| > k\}$ and take $G_k(u_p)$ $(k>0)$  as test function in \eqref{Pp2}, yielding to   
	\begin{equation}\label{stima1BIS}
	\int_{A_k} \phi(|\nabla u_p|)|\nabla u_p|^{p+1}  dx = \int_\Omega fG_k(u_p) dx\leq \mathcal{S}_1 \|f\|_{L^{N}(\Omega)} \int_{A_k} |\nabla u_p| dx,
	\end{equation}
	after applications of the H\"older and the Sobolev inequalities.
	
	Observe that  by \eqref{phi3} one has that, for a fixed  $0<\varepsilon<1$ there exists $ {h}$ such that 
	\begin{equation}\label{controllobassoBIS}
		\phi(|\nabla u_p|)|\nabla u_p|^{p+1}\geq (1-\varepsilon) |\nabla u_p|^p \ge (1-\varepsilon) |\nabla u_p| - (1-\varepsilon)\frac{p-1}{p} \,\ \ \text{on}\ \ \{|\nabla u_p|>{ {h}}\}, 
	\end{equation}
	where we also used the Young inequality. Then, gathering \eqref{controllobassoBIS} with \eqref{stima1BIS}, one gets
	\begin{equation}\label{stima2BIS}
		\begin{aligned}
		(1-\varepsilon)\int_{A_k \cap \{|\nabla u_p|>h\}} |\nabla u_p| dx & \leq \mathcal{S}_1 \|f\|_{L^{N}(\Omega)} \int_{A_k \cap \{|\nabla u_p|>h\}} |\nabla u_p| dx + \mathcal{S}_1 \|f\|_{L^{N}(\Omega)}h|A_k|
		\\
		& + (1-\varepsilon)\frac{p-1}{p}|A_k|.
		\end{aligned}	
	\end{equation}
	
	Now we fix $\overline{\varepsilon}$  such that 
	$$
	1- \overline{\varepsilon} -\mathcal{S}_1 \|f\|_{L^{N}(\Omega)}>0,
	$$
	and then we can fix $\overline{h}:=\overline{h}(\overline{\varepsilon})$ such that \eqref{controllobassoBIS} holds. 
	Then from \eqref{stima2BIS} one obtains 
	\begin{equation}\label{hoveline}
		\begin{aligned}
			(1- \overline{\varepsilon} -\mathcal{S}_1 \|f\|_{L^{N}(\Omega)})\int_{A_k\cap \{|\nabla u_p|>\overline{h}\}}|\nabla u_p|  dx \leq \left(\mathcal{S}_1 \|f\|_{L^{N}(\Omega)}\overline{h} + \frac{p-1}{p}\right)|A_k|,
		\end{aligned}
	\end{equation}	
	where the right-hand of \eqref{hoveline} can be estimated independently of $p$ as $\frac{p-1}{p}< 1$. 
	
	Therefore
	\begin{eqnarray*}
		(1- \overline{\varepsilon} -\mathcal{S}_1 \|f\|_{L^{N}(\Omega)})\int_{\Omega}|\nabla G_k( u_p)|  dx  &\leq  &(1- \overline{\varepsilon} -\mathcal{S}_1 \|f\|_{L^{N}(\Omega)}) \left(\int_{A_k\cap \{|\nabla u_p|>\overline{h}\}}|\nabla u_p|dx +|A_k|\overline{h} \right)\\ 
		&\stackrel{\eqref{hoveline}}{\leq} & (\mathcal{S}_1 \|f\|_{L^{N}(\Omega)}\overline{h}+1 + \overline{h})|A_k| \\
		&  \stackrel{\eqref{hf}}{\leq} & (2\overline{h}+1)|A_k|, 
	\end{eqnarray*}
	which implies that
	\begin{equation}\label{stima3}
		\int_\Omega |\nabla G_k(u_p)| dx  \le \frac{(2\overline{h}+1)|A_k|}{1-  {\overline{\varepsilon}} -\mathcal{S}_1 \|f\|_{L^{N}(\Omega)}}.
	\end{equation}
	Now one can apply to \eqref{stima3} both the Sobolev and  the H\"older inequalities, deducing that 
	\begin{equation*}
		\int_\Omega |G_k(u_p)| dx \le \frac{(2\overline{h}+1)\mathcal{S}_1|A_k|^{1+\frac{1}{N}}}{1- \overline{\varepsilon} -\mathcal{S}_1 \|f\|_{L^{N}(\Omega)}}.
	\end{equation*}
	In particular,  for any $\overline{k}>k>0$, one has 
	\begin{equation}\label{stima4}
		|A_{\overline{k}}| \le \frac{(2\overline{h}+1)\mathcal{S}_1 |A_k|^{1+\frac{1}{N}}}{(\overline{k}-k)(1- \overline{\varepsilon}-\mathcal{S}_1 \|f\|_{L^{N}(\Omega)})}
	\end{equation}	
	which allows to apply a standard Stampacchia argument (see \cite{Stampacchia}) and so deduce
	\begin{equation}\label{stimalimitatezza}
		\|u_p\|_{L^\infty(\Omega)}\le M,
	\end{equation}
	where $M>0$ does not depend on $p$ as the right-hand of \eqref{stima4} does not.
	
	As for the estimate \eqref{3222}, we take $u_p$ as test function in \eqref{Pp2} and, using Lemma \ref{coerc} one  gets
	\[\alpha \int_{\Omega}|\nabla u_p|^pdx-B_p|\Omega|\le \int_{\Omega}\phi(|\nabla u_p|)|\nabla u_p|^{p+1}dx=\int_{\Omega}fu_p\, dx\stackrel{\eqref{stimalimitatezza}}{\le} M\| f\|_{L^1(\Omega)}\,.\]
	 Since $B_p$ is bounded (recall we may assume $p$ small without loosing generality), estimate \eqref{3222} follows. Finally   the $BV$-estimate is now a consequence of Young's inequality.
\end{proof}

Next result concerns the existence of the vector field $\z \in L^\infty(\Omega,\mathbb{R}^N)$.  
\begin{corollary}
\label{corUZ2}
Let $\phi$ satisfy \eqref{phi1}, \eqref{phi2} and \eqref{phi3} and let $f \in L^N(\Omega)$ be such that \eqref{hf} is in force. Let $u_p$ be a solution to \eqref{Pp2}. Then there exists $u \in BV(\Omega)\cap L^\infty(\Omega)$ such that, up to subsequences, $u_p \to u$ in $L^q(\Omega)$ for any $q<\infty$ as $p\to 1^+$.  Moreover there exists $\z\in L^\infty(\Omega, \mathbb{R}^N)$ with $\|\z\|_\infty\le 1$ such that 
\begin{equation}\label{toz3}
	\phi(|\nabla u_p|)|\nabla u_p|^{p-1}\nabla u_p {\rightharpoonup} \, \z \quad \mbox{in $L^q(\Omega, \mathbb{R}^N)$ for every $1\le q<\infty$}.
\end{equation}
\end{corollary}
\begin{proof}
	It follows from Lemma \ref{lemma.boundedness2} that $u_p$ is bounded in $BV(\Omega)\cap L^\infty(\Omega)$ with respect to $p$. Then standard compactness arguments apply, yielding a function $u\in BV(\Omega)\cap L^\infty(\Omega)$ such that, up to subsequences, $u_p \to u \mbox{ in $L^q(\Omega)$}$ for any $1 \le q <\infty$. 	
	
	Now let $q<\frac{p}{p-1}$ and observe that it follows from an application of the H\"older inequality and from \eqref{phi3} that 
	 \begin{equation}\label{zBIS}
	 \begin{aligned}	
	\int_\Omega \Big|\phi(|\nabla u_p|)|\nabla u_p|^{p-1}\nabla u_p\Big|^q dx \le \int_\Omega |\nabla u_p|^{q(p-1)}dx 
	&\le  \left(\int_\Omega |\nabla u_p|^p dx \right)^{\frac{q(p-1)}{p}}|\Omega|^{1-\frac{q(p-1)}{p}}
	\\
	&\stackrel{\eqref{3222}}{\leq}  C^{\frac{q(p-1)}{p}} |\Omega|^{1-\frac{q(p-1)}{p}}.
	\end{aligned}		
	\end{equation}
	Estimate \eqref{zBIS} gives, after a standard diagonal argument, that  there exists a unique vector field $\z$ such that \eqref{toz3} holds.
	Moreover, reasoning as in \cite{MRST},  by the lower semicontinuity with respect to $p$ and by taking $q\to\infty$ into \eqref{zBIS}, one obtains that $\z \in L^\infty(\Omega,\mathbb{R}^N)$ and  $\|\z\|_{\infty}\le 1$. This concludes the proof.
\end{proof}

\subsection{Passage to the limit}

In this section we prove Theorem \ref{theorem1} by gathering a series of lemmas.  The first one addresses the identification of the vector field $\z$.

\begin{lemma}
\label{lemmaTeo3}
	Let $\phi$ satisfy \eqref{phi1}, \eqref{phi2} and \eqref{phi3} and let $f \in L^N(\Omega)$ be such that \eqref{hf} is in force. Then  $u\in BV(\Omega)\cap L^\infty(\Omega)$ and $\z \in L^\infty(\Omega, \mathbb{R}^N)$ found in Corollary \ref{corUZ2} satisfy:
	\begin{equation}\label{lemTeo3_divz}
	-\operatorname{div}\z = f(x) \quad \mbox{in $\mathcal{D}'(\Omega)$},
\end{equation}
and 
\begin{equation}\label{lemTeo4_pairing}
	(\z,Du) = \Phi(|Du|) +\Phi^*(|\z|)\mathcal{L}^N \quad \mbox{in $\mathcal{M}(\Omega)$}.
\end{equation}

\end{lemma}
\begin{proof}
Let $u_p$ be a solution to \eqref{Pp2}. It follows immediately from Corollary \ref{corUZ2} that one can easily pass to the limit into the weak formulation of \eqref{Pp2} in order to deduce that \eqref{lemTeo3_divz} holds. 

We now turn to show \eqref{lemTeo4_pairing}. 
Let $\varphi\in C^1_c(\Omega)$ be non-negative and take $u_p\varphi$ as test function in \eqref{Pp2}. Then
\begin{equation}\label{pairing2}
	 \int_\Omega \varphi \phi(|\nabla u_p|)|\nabla u_p|^{p+1} dx 
	+ \int_\Omega u_p\phi(|\nabla u_p|)|\nabla u_p|^{p-1}\nabla u_p\nabla \varphi\, dx = \int_\Omega fu_p\varphi\, dx.
\end{equation} 
One can pass to the limit on the right-hand of \eqref{pairing2} as, by Corollary \ref{corUZ2}, $u_p \to  u$ in $L^q(\Omega)$ for $1\le q<\infty$ and $f\varphi\in L^N(\Omega)$. We also may pass to the limit in the second term on the left-hand thanks to Corollary \ref{corUZ2}. Our next concern is to show that we may also pass to the limit in the first term.

The first step is to write that term in a suitable way. To this aim, set $\z_p=\phi(|\nabla u_p|)\nabla u_p$. If $\nabla u_p\ne0$, then (by \eqref{YoungFencheleq1})
\[\z_p\cdot\nabla u_p=\Phi(|\nabla u_p|)+\Phi^*(|\z_p|),\]
so that 
\[\z_p\cdot\nabla u_p |\nabla u_p|^{p-1}=\Phi(|\nabla u_p|)|\nabla u_p|^{p-1}+\Phi^*(|\z_p|)|\nabla u_p|^{p-1}.\]
Observe that this identity is also satisfied when $\nabla u_p=0$.
As a consequence, it holds
\begin{equation}
\begin{aligned}\label{Ident_Pairing}
\int_\Omega \varphi\phi(|\nabla u_p|)|\nabla u_p|^{p+1}dx &=\int_\Omega \varphi \z_p\cdot \nabla u_p |\nabla u_p|^{p-1}dx\\
&=\int_\Omega \varphi \Phi(|\nabla u_p|)|\nabla u_p|^{p-1}dx+\int_\Omega \varphi\Phi^*(|\z_p|) |\nabla u_p|^{p-1}dx
\\
&=I_1+I_2.
\end{aligned}
\end{equation}
We separately analyze the last two terms in order to let $p\to 1^+$. For $I_1$, we fix $k>0$ and perform the following calculations
\begin{equation*}
\begin{aligned}
\int_\Omega \varphi \Phi(|\nabla u_p|)|\nabla u_p|^{p-1}dx &\ge 
\int_\Omega \varphi G_k\big(\Phi(|\nabla u_p|)\big)|\nabla u_p|^{p-1}dx\\
&\ge 
\big(\Phi^{-1}(k)\big)^{p-1}\int_\Omega \varphi (G_k\big(\Phi(|\nabla u_p|)\big))\, dx,
\end{aligned}
\end{equation*}
since $\Phi^{-1}(k) = \inf \{s \geq 0: \, \Phi(s) > k\}$.
Having in mind that $G_k\circ\Phi$ is a convex real function which induces a lower semicontinuous operator (see \cite[Fact 3.2]{anzellotti}), we obtain
\begin{equation}\label{G_k-measure}
\liminf_{p\to 1^+}\int_\Omega \varphi \Phi(|\nabla u_p|)|\nabla u_p|^{p-1}dx\ge 
\liminf_{p\to 1^+} \int_\Omega \varphi G_k\big(\Phi(|\nabla u_p|)\big)\, dx
\ge \int_\Omega \varphi G_k\big(\Phi(|Du|)\big).
\end{equation}
To identify the measure $G_k\big(\Phi(|Du|)\big)$, observe that, for every $\xi\in\R^N$,
\[\lim_{t\to\infty}\frac{G_k\big(\Phi(t|\xi|)\big)}t=
\lim_{t\to\infty}\left(\frac{\Phi(t|\xi|)}t-\frac kt\right)^+=|\xi|\,.
\]
Thus, 
\[\int_\Omega \varphi G_k\big(\Phi(|Du|)\big)=\int_\Omega \varphi G_k\big(\Phi(|D^au|)\big)\, dx+\int_\Omega \varphi |D^su|\,,\]
so that, letting $k\to 0^+$, \eqref{G_k-measure} becomes
\begin{equation}\label{G_k-measure1}
\liminf_{p\to 1^+}\int_\Omega \varphi \Phi(|\nabla u_p|)|\nabla u_p|^{p-1}dx
\ge \int_\Omega \varphi \Phi(|Du|)\,.
\end{equation}

On the other hand, we deal with term $I_2$ by considering $\alpha=\lim_{s\to0^+}\phi(s)s$. If $\alpha=1$, then the associated operator is the 1--Laplacian, for which $\Phi^*(s)=0$ for all $0\le s\le1$ and so there is nothing to check. Hence, assume that $0\le\alpha<1$.  
 
In particular we want to show that
\begin{equation}
	\liminf_{p \to 1^+} \int_{\Omega} \varphi \Phi^*(|\z_p|) |\nabla u_p|^{p-1} \, dx \ge \int_{\Omega} \varphi \Phi^*(|\z|) dx.
	\label{eq:goal}
\end{equation}
Now observe that it follows from the definition of $\Phi^*$ that for every $\xi \in \mathbb{R}^N$ it holds
\begin{equation}\label{eqdef}
	\Phi^*(|\z_p|) \ge \xi \cdot \z_p - \Phi(|\xi|).
\end{equation}
Then, let us set 
\begin{equation*}
	K_p(x) = \frac{1}{p} + \frac{p-1}{p} |\nabla u_p(x)|^p,
\end{equation*} 
and let us multiply \eqref{eqdef} by $|\nabla u_p(x)|^{p-1}$, obtaining
\begin{equation}\label{eqgoal0}
	\Phi^*(|\z_p|) |\nabla u_p|^{p-1} \ge \xi \cdot \z_p |\nabla u_p|^{p-1} - \Phi(|\xi|) |\nabla u_p|^{p-1} \ge \xi \cdot \z_p |\nabla u_p|^{p-1} - K_p\Phi(|\xi|), 
\end{equation}
where the last inequality follows from an application of the Young inequality. 
Then, as it holds
\begin{equation*}
	K_p \sup_{\xi \in \mathbb{R}^N} \left\{ \xi \cdot \frac{\z_p |\nabla u_p|^{p-1}}{K_p} - \Phi(|\xi|) \right\} = K_p \Phi^* \left( \frac{|\z_p||\nabla u_p|^{p-1}}{K_p} \right)
\end{equation*}
from \eqref{eqgoal0} one gains
\begin{equation*}
	\Phi^*(|\z_p|) |\nabla u_p|^{p-1} \ge K_p \Phi^* \left( \frac{|\z_p||\nabla u_p|^{p-1}}{K_p} \right) \quad \text{a.e. in } \Omega,
\end{equation*}
that is 
\begin{equation*}
	\int_{\Omega} \varphi \Phi^*(|\z_p|) |\nabla u_p|^{p-1} \, dx \ge \int_{\Omega} \varphi K_p \Phi^* \left( \frac{|\z_p||\nabla u_p|^{p-1}}{K_p} \right) \, dx.
\end{equation*}
Now, as it is immediate to show that $K_p$ converges to $1$ in $L^1(\Omega)$ and as $\z_p|\nabla u_p|^{p-1} \rightharpoonup \z$, it follows from Corollary $3.9$ in \cite{Br} that it holds
\begin{equation}
	\begin{aligned}
		\liminf_{p \to 1^+} \int_{\Omega} \varphi K_p \Phi^* \left( \frac{|\z_p||\nabla u_p|^{p-1}}{K_p} \right) \, dx &\ge \liminf_{p \to 1^+} \frac{1}{p}\int_{\Omega} \varphi \Phi^* \left( \frac{|\z_p||\nabla u_p|^{p-1}}{K_p} \right) \, dx
		\\
		&\ge \int_{\Omega} \varphi \Phi^* \left(|\z| \right) \, dx,
	\end{aligned}	
\end{equation}
which gives \eqref{eq:goal}.

%We remark that the generalized inverse defined in \eqref{inverse}, satisfies $\psi(s)=0$ for all $0\le s\le\alpha$, wherewith $\Phi^*(s)=0$ for all $0\le s\le\alpha$. Thus, $t>\alpha$ if and only if $\Phi^*(t)>0$. Now for a fixed $k>0$, notice that
%\begin{equation}\label{inverse2}
%\Phi^*(|\z_p|)>k \Longleftrightarrow \phi(|\nabla u_p|)|\nabla u_p|>(\Phi^*)^{-1}(k).
%\end{equation}
%Observe that $\phi(s)s>t$ implies $s>\psi(t)$. As a consequence, denoting $h=\psi\left((\Phi^*)^{-1}(k)\right)$, it follows from \eqref{inverse2} that $|\nabla u_p|>h$. Moreover, on account of $(\Phi^*)^{-1}(k)>\alpha$, it yields $h>0$. Therefore,
%\[\int_\Omega\varphi\Phi^*(|\z_p|)|\nabla u_p|^{p-1}dx\ge
%\int_\Omega\varphi G_k\big(\Phi^*(|\z_p|)\big)|\nabla u_p|^{p-1}dx\ge \psi\left(\big(\Phi^*\big)^{-1}(k)\right)^{p-1}\int_\Omega\varphi G_k\big(\Phi^*(|\z_p|)\big)\, dx.\]
%Letting $p\to 1^+$, we deduce from \cite[Corollary 3.9]{Br} that
%$$
%\begin{aligned}
%\liminf_{p\to1^+}\int_\Omega\varphi\Phi^*(|\z_p|)|\nabla u_p|^{p-1}dx&\ge
%\liminf_{p\to1^+}h^{p-1}\int_\Omega\varphi G_k\big(\Phi^*(|\z_p|)\big)\, dx
%\\
%&\ge \liminf_{p\to1^+}\int_\Omega\varphi G_k\big(\Phi^*(|\z_p|)\big)\, dx
%\ge \int_\Omega\varphi G_k\big(\Phi^*(|\z|)\big)\, dx.
%\end{aligned}
%$$
%Since this inequality holds for every $k>0$, we obtain that
%\begin{equation}\label{G_k-measure2}
%\liminf_{p\to1^+}\int_\Omega\varphi\Phi^*(|\z_p|)|\nabla u_p|^{p-1}dx
%\ge \int_\Omega\varphi \Phi^*(|\z|)\, dx.
%\end{equation}

Now, going back to \eqref{Ident_Pairing} and taking into account \eqref{G_k-measure1} and \eqref{eq:goal}, we conclude that 
\begin{equation}\label{anInequality}
\liminf_{p\to 1^+}\int_\Omega \varphi\phi(|\nabla u_p|)|\nabla u_p|^{p+1}dx\ge 
\int_\Omega \varphi \Phi(|Du|)+\int_\Omega\varphi \Phi^*(|\z|)\, dx.
\end{equation}
Inequality  \eqref{anInequality} is the last step to let $p\to 1^+$ in \eqref{pairing2}, so that it yields
\begin{equation*}
	 \int_\Omega \varphi \Phi(|Du|) +\int_\Omega \varphi \Phi^*(|\z|) 
	+ \int_\Omega u\z\cdot \nabla \varphi\, dx \le \int_\Omega fu\varphi\, dx.
\end{equation*} 
Having in mind \eqref{lemTeo3_divz}, the last inequality can be written as
\begin{equation*}
	 \int_\Omega \varphi \Phi(|Du|) +\int_\Omega \varphi \Phi^*(|\z|) 
	\le - \int_\Omega u\z\cdot \nabla \varphi\, dx - \int_\Omega u\varphi\operatorname{div}\z\, dx=\int_\Omega\varphi (\z, Du)\,,
\end{equation*} 
wherewith $\Phi(|Du|) +\Phi^*(|\z|) 	\le  (\z, Du)$ holds as measures.
As the reverse inequality follows from \eqref{YoungFencheleq2},  the proof is complete.
\end{proof}

Now we focus on the boundary condition.

\begin{lemma}
	\label{lemma.boundaryBIS}
	Let $\phi$ satisfy \eqref{phi1}, \eqref{phi2} and \eqref{phi3} and let $f \in L^N(\Omega)$ be such that \eqref{hf} is in force. Then the function $u$ and the vector field $\z$ defined in Corollary \ref{corUZ2} satisfy
	$$
		[\z,\nu]\in \sgn(-u) \ \text{on } \partial\Omega.
	$$
\end{lemma}
\begin{proof}
	Let $u_p$ be a solution to \eqref{Pp2} and let us take $u_p$ as a test function in the weak formulation of \eqref{Pp2}, yielding to
	\begin{equation}\label{bordo1BIS}
		\int_\Omega \phi(|\nabla u_p|)|\nabla u_p|^{p+1} dx  =  \int_{\Omega}fu_p dx.
	\end{equation}
	We aim to take the $\liminf$ as $p\to 1^+$ into the previous equality.
	For the right-hand of \eqref{bordo1BIS} we can simply pass to the limit as $p\to 1^+$ since $f\in L^N(\Omega)$ and $u_p$ converges to $u$ in $L^q(\Omega)$ for every $1 \le q < \infty$ (see Corollary \ref{corUZ2}). For the left-hand side of \eqref{bordo1BIS} we can reason similarly to what was done to obtain \eqref{G_k-measure1} and \eqref{eq:goal} (just replacing \eqref{eq.lsc1} by \eqref{eq.lsc2}), yielding to
\begin{equation}\label{anInequality2}
\liminf_{p\to 1^+}\int_\Omega \phi(|\nabla u_p|)|\nabla u_p|^{p+1}dx\ge 
\int_\Omega \Phi(|Du|) + \int_{\partial \Omega}|u| d\mathcal{H}^{N-1} + \int_\Omega \Phi^*(|\z|)\, dx.
\end{equation}	
Hence, from \eqref{bordo1BIS} and \eqref{anInequality2}, it follows that	
	\begin{equation*}\label{bordo2BIS}
		\begin{aligned}
			\int_\Omega \Phi(|Du|) + \int_{\partial \Omega} |u| d\mathcal{H}^{N-1} + \int_\Omega \Phi^*(|\z|)dx &\le  \int_{\Omega}fu dx \stackrel{\eqref{lemTeo3_divz}}{=}-\int_\Omega u \operatorname{div}\z dx
			\\
			& \stackrel{\eqref{green}}{=} \int_\Omega (\z,Du) - \int_{\partial\Omega} u [\z,\nu] d\mathcal{H}^{N-1},
		\end{aligned}
	\end{equation*}
	which, thanks to \eqref{lemTeo4_pairing}, takes to 
	$$
	\int_{\partial\Omega} \left(|u|+u[\z,\nu] \right)\,d\mathcal{H}^{N-1}\leq0.
	$$
	This concludes the proof since $|[\z,\nu]|\le 1$.
\end{proof}

Hence, gathering together Lemmas \ref{lemmaTeo3} and \ref{lemma.boundaryBIS}, we complete the proof of Theorem \ref{theorem1}.

\section{On the identification of the vector field $\z$}
\label{s5}

For solutions in the sense of Definition \ref{def} the vector field $\z$ is not, in general, uniquely determined. This is a well known fact, for instance, for equations involving the 1-Laplace operator (i.e. $\phi(s)s=1$). In other cases, as for the non-parametric minimal surface equation $\left(\phi(s)s=\frac{s}{\sqrt{1+s^2}}\right)$, the vector field can be uniquely determined as a function of the regular part of $D u$. 

\medskip

In this section we want to briefly address this phenomenon; in particular we are interested in understanding whether the vector field  $\z$ can be derived as a function of $D^au$. As we will see this is strictly related to the structure of the function $\phi$.     

\medskip

In the entire paper, $\phi$ is a function defined on $[0,\infty)$ which, as \eqref{phi1} is in force, is continuous in $(0,\infty)$. 
For the next result, we strengthen the hypotheses on $\phi$ requiring
\begin{equation} \label{typeI}
	\lim_{s\to0^+}\phi(s)s=0,
\end{equation}	
providing a function $\Phi$ which is $C^1(\mathbb{R})$. As we rely on \eqref{YoungFencheleq1}, this is sufficient to have an explicit representation for $\z$ as stated in the next proposition.
\begin{proposition}\label{propunique}
	Assume that $\phi$ satisfies \eqref{phi1}, \eqref{phi2}, \eqref{phi3} and \eqref{typeI}.
	Then any solution to \eqref{P1} in the sense of Definition \ref{def} is such that $\z$ has the form
	\begin{equation}\label{uniquez}
		\z=\phi(|D^a u|)D^a u.
	\end{equation} 
\end{proposition}
\begin{proof}
	Let $u$ be a solution to \eqref{P1} in the sense of Definition \ref{def} and let $\z$ be the associated vector field. First recall that it follows from \eqref{eqweakform2} that (recall \eqref{YoungFenchelIdweak1}):
	\begin{equation*} \label{eqassoluta}
		\z \cdot D^a u = \Phi(|D^a u|) + \Phi^{*}(|\z|) \text{ in }\Omega,
	\end{equation*}
	which is the equality case in the Fenchel-Young inequality, namely \eqref{YoungFencheleq1}. 
	Now it is sufficient to observe that, thanks to \eqref{phi1} and to \eqref{typeI}, $\Phi \in C^1(\mathbb{R})$. Then, recalling once again the discussion in Section \ref{FYI}, one gains that 
	\begin{equation*}
		\z = \phi(|D^a u|)D^a u,
	\end{equation*}
	almost everywhere in $\Omega$. This completes the proof.
\end{proof}

Let us spend a few words on Proposition \ref{propunique}.

\begin{remark}
	\rm Proposition \ref{propunique} gives a formula to represent $\z$ from the regular part of the gradient $D^a u$. 
	A natural question is whether one can derive $D^a u$ from $\z$. This is clearly possible under the stronger assumption that $\phi(s)s$ is strictly increasing which, in this case, provides a complete knowledge on $\z$ and of his role in Definition \ref{def}.
\end{remark}

\begin{remark}
		\rm As for the Proposition \ref{propunique}, we have in mind essentially two types of functions whose graphs are given in Figure \ref{f1} and \ref{f2} below respectively.   
		\begin{figure}[h]
			\centering
			\begin{tikzpicture}
				\draw [->,thin] (-0.2,0) -- (9,0);
				\draw [->,thin] (0,-0.2) -- (0,1.5);
				\draw [thick, blue] (0.05,0) -- (1,0);
				\draw [thick, domain=1:3, smooth, variable=\x, blue] plot({\x}, {(-1+\x)/(1+\x)});
				\draw [thick, blue] (3,0.5) -- (6,0.5);
				\draw [thick, domain=6:9, smooth, variable=\x, blue] plot({\x}, {(-4+\x)/(-2+\x)});
				\fill [blue] (0,0) circle (1.5pt);
				\node [left] at (0,1.3) { $\phi(s)s$};
				\node [below] at (8.8,0) { $s$};
			\end{tikzpicture}
			\caption{Functions $\phi(s)s$ satisfying \eqref{typeI}}
			\label{f1}
		\end{figure}
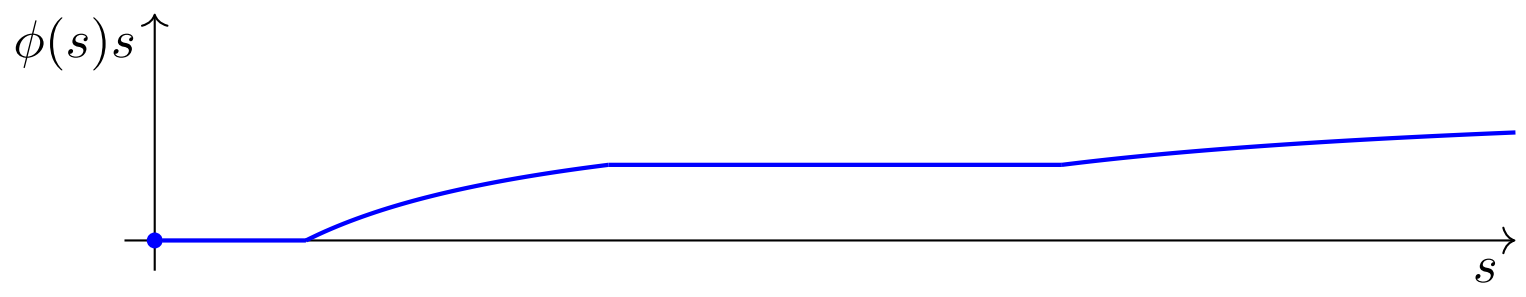 
		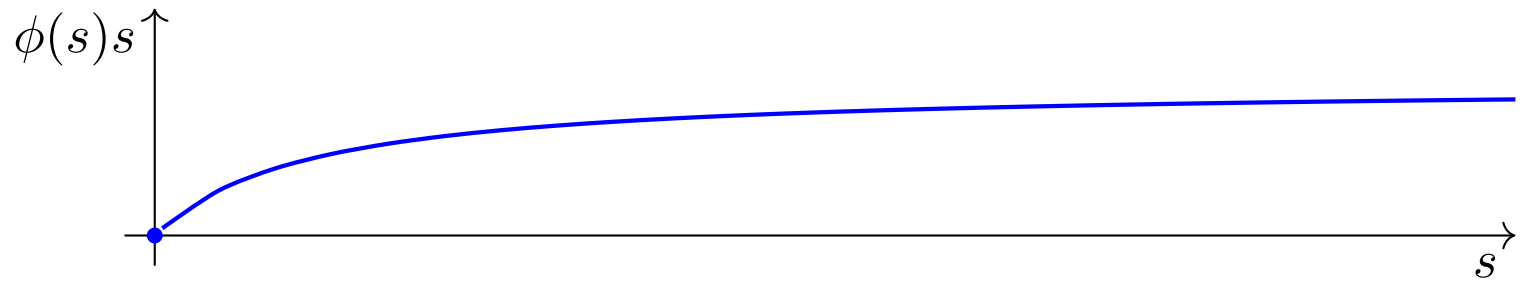
\begin{figure}[h]
			\centering
			\begin{tikzpicture}
				\draw [->,thin] (-0.2,0) -- (9,0);
				\draw [->,thin] (0,-0.2) -- (0,1.5);
				\draw [thick, domain=0.05:9, smooth, variable=\x, blue] plot({\x}, {\x/(1+\x)});
				\fill [blue] (0,0) circle (1.5pt);
				\node [left] at (0,1.3) { $\phi(s)s$};
				\node [below] at (8.8,0) { $s$};
			\end{tikzpicture}
			\caption{Functions $\phi(s)s$ strictly increasing satisfying \eqref{typeI}}
			\label{f2}
		\end{figure}
		When $\phi(s)s$ satisfies \eqref{typeI}, then Proposition \ref{propunique} holds and the formula \eqref{uniquez} completely characterizes $\z$ and $D^a u$ if $\phi(s)s$ is strictly increasing. 
	
		\medskip
		
		For functions of the above type $\Phi^*(t)>0$ for $t \ne 0$, while either $\Phi(s)>0$ for $s\ne 0$ or $\Phi(s)=0$ for $s$ near $0$ according to, respectively, a positive or null $\phi(s)s$ for s near $0$. The mean curvature operator is the prototypical example of this class;  here $\phi(s)s=\frac{s}{\sqrt{1+s^2}}$ and one has $$(\z,Du)=\sqrt{1+|Du|^2} - \sqrt{1-|\z|^2} \mathcal{L}^N.$$ 
		In this case $\Phi(s)$ is increasing and $C^1(\mathbb{R})$ and $\z$ is uniquely determined by
		$$\z= \frac{D^a u}{\sqrt{1+|D^a u|^2}}.$$ 
		
		\medskip
		
		On the other side, we also deal with functions $\phi(s)s$ which are not continuous up to the origin, namely $\lim_{s\to 0^+} \phi(s)s > 0$. In this case, $\Phi\in C^1((0,\infty))$, and the relationship between $\z$ and $Du$ cannot be described by a single formula, i.e. Proposition \ref{propunique} does not apply. In this context, a unique representation for $\z$ cannot be derived; however, a suitable notion of subdifferential $\partial\Phi$ can be set so that one can still give a representation formula for $\z$. Just to give an idea, in the case of the $1$-Laplacian where $\phi(s)s=1$ for $s>0$, which is the prototype for this class of functions, and reasoning similarly to the proof of Proposition \ref{propunique}, one has that if $D^a u \neq 0$, then $\z$ is uniquely given by $\z = D^a u / |D^a u|$. Otherwise, if $D^a u = 0$, the relation reduces to the inequality $|\z| \le 1$ as $[-1,1] = \partial \Phi(0)$.
		
		\medskip
		
		The functions just described are the ones with graph given by 
		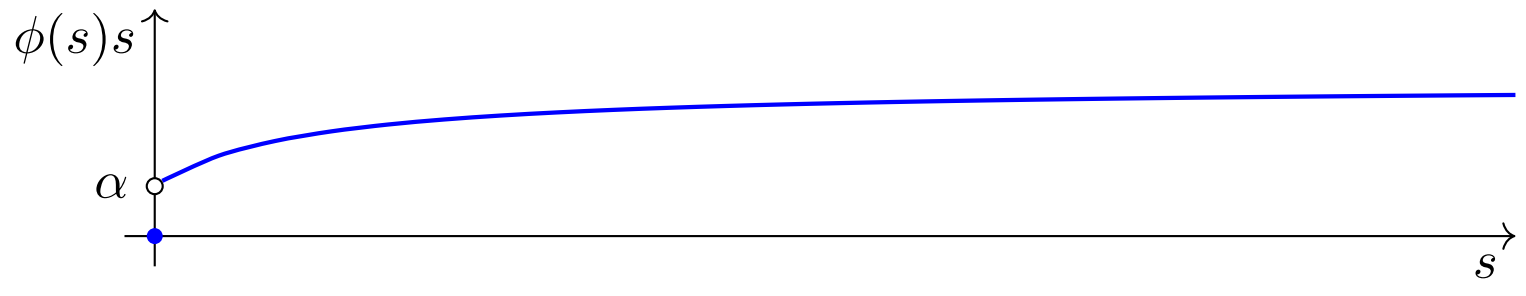
\begin{figure}[h]
			\centering
			\begin{tikzpicture}
				\draw [->,thin] (-0.2,0) -- (9,0);
				\draw [->,thin] (0,0.38) -- (0,1.5);
				\draw [thin] (0,-0.2) -- (0,0.28);
				\draw [thick, domain=0.05:9, smooth, variable=\x, blue] plot({\x}, {(1+3*\x)/(3+3*\x)});
				\fill [blue] (0,0) circle (1.5pt);
				\draw (0,0.33) circle (1.5pt);
				\node [left] at (-0.05,0.33) {$\alpha$};
				\node [left] at (0,1.3) { $\phi(s)s$};
				\node [below] at (8.8,0) {$s$};
			\end{tikzpicture}
			\caption{Functions $\phi(s)s$ not satisfying \eqref{typeI}}
			%\label{,,}
		\end{figure}
		\\ Here $\Phi^*(t)=0$ for $|t|\le \alpha$, while $\Phi(s)>0$ for all $s\ne0$. In the general case, where $\alpha = \lim_{t \to 0^+}\phi(s)s$, if $D^a u \neq 0$, then $\z$ is uniquely given by $\z = \phi(|D^a u|)D^a u$. Otherwise, if $D^a u = 0$, the relation reduces to the inequality $|\z| \le \alpha$ as $[-\alpha,\alpha] = \partial \Phi(0)$.
\end{remark}	

As a consequence of Proposition \ref{propunique} and the above remarks, one has the following: 
\begin{corollary}
Assume that $\phi$ satisfies \eqref{phi1}, \eqref{phi2}, \eqref{phi3} and $\phi(s)s$ is increasing. Then $D^au$ is determined by $\z$.
\end{corollary}
\begin{proof} Denote  by 
$\alpha=\lim_{s\to0^+} \phi(s)s$ and let $\psi (s)$ the inverse function of $\phi(s)s$. 
One easily  realizes that 
$|D^a u|=\psi(|\z|)$ if $|\z|>\alpha$ and $|D^a u|=0$ if $|\z|\le \alpha$.
\end{proof}

\section{The case of the Lorentz space $ L^{N,\infty}(\Omega)$} 

In this section we briefly deal with the case of a datum $f$ belonging to the Lorentz space $L^{N,\infty}(\Omega)$. 

Let us briefly recall (see  \cite{PKF} for further details) that the decreasing rearrangement of a measurable function $v$, is denoted by $v^*$ and it is  defined as  $$
	v^{*}(s)=\sup\{t>0: |\{|v|>t\}|>s\}.
	$$
	For fixed $p>1$ and  $ 1 \leq q \leq \infty $,   let  $v$ be a measurable function.   We say that $v$ belongs to the $(p,q)$-Lorentz space if the following quantity is finite
	\begin{equation*}\label{deflorentz}
		\|v\|_{L^{p,q}(\Omega)}=  \begin{cases}
			  \left(\int_{0}^{\infty}[t^{1/p}v^{*}(t)]^q\frac{\,dt}{t}\right)^{\frac{1}{q}} & 1<p<\infty,\,1\leq q<\infty, \\
			  \sup_{t>0} t^{\frac{1}{p}}v^{*}(t) & 1<p<\infty,\,q=\infty;
		\end{cases}
	\end{equation*}
	we also recall that  $L^{p,p}(\Omega)=L^{p}(\Omega)$. 
	\\Moreover, considering \( 1 < q < p < r < \infty \), the following inclusions are easy to check:
	$$
	L^r(\Omega)\subset L^{p,1}(\Omega)\subset L^{p,q}(\Omega)\subset  L^{p}(\Omega)\subset L^{p,r}(\Omega)\subset L^{p,\infty}(\Omega)\subset L^{q}(\Omega).
	$$
 To our purpose we just mention that a Sobolev inequality is available, namely
\begin{equation}\label{soblorentz}
	\|v\|_{L^{1^{\ast},1}(\Omega)}\leq \tilde{\mathcal S}_{1} \|v\|_{W^{1,1}_{0}(\Omega)}, \ \ \forall v\in W^{1,1}_{0}(\Omega),
\end{equation}
where
$$ \tilde{\mathcal S}_{1}= [(N-1)\omega_{N}^{\frac{1}{N}}]^{-1}. \ \ $$

%First of all we highlight that Theorem \ref{theorem1} still holds in this case. 

Before stating the existence result, let us stress that solutions to \eqref{P1} with $f\in L^{N,\infty}(\Omega)$ are meant in the sense of Definition \ref{def}, where the requirement $\z\in X_N(\Omega)$ is replaced by $\z\in X_{N,\infty}(\Omega):=\{\z\in L^\infty(\Omega,\mathbb{R}^N):\ \operatorname{div}\z\in L^{N,\infty}(\Omega)\}$. We also point out that the Anzellotti theory, recalled in Section \ref{sec:anz}, can be shown to be still available in this framework; in particular Green formula \eqref{green} continues to hold when $v\in BV(\Omega)$ and $\z \in X_{N,\infty}(\Omega)$.  
Hence, with minor modifications in the proof of Theorem \ref{theorem1}, one can prove the following existence result: 

\begin{theorem}\label{theorem_marc}
	Let $\phi$ satisfy \eqref{phi1}, \eqref{phi2} and \eqref{phi3} and let $f \in L^{N,\infty}(\Omega)$ be such that
	\begin{equation}\label{hf2}
		\|f\|_{L^{N,\infty}(\Omega)}<\tilde{\mathcal S}_{1}^{-1}.
	\end{equation} 
	Then there exists  a bounded solution of \eqref{P1}.
\end{theorem}
As a further property we show the following degeneration phenomenon:
\begin{theorem}\label{thm:all-solutions-trivial}
Assume that $\phi$ satisfies \eqref{phi1}, \eqref{phi2}, \eqref{phi3} and set $\alpha=\lim_{s\to 0^+}\phi(s)s>0$.
If $\| f\|_{L^{N,\infty}(\Omega)}\mathcal <\alpha\tilde{S}_1^{-1}$, then any solution to \eqref{P1} vanishes identically.
\end{theorem}

\begin{proof}
	Let $u$ be any solution of \eqref{P1} and let
	$\z\in L^\infty(\Omega,\mathbb{R}^N)$, $\|\z\|_\infty\le1$, be an associated vector
	field.
	Choosing $v=u$ in \eqref{green} and using \eqref{eqweakform1}, we get
	\[
	-\int_\Omega f u\,dx+\int_\Omega (\z,Du)=\int_{\partial\Omega}u\,[\z,\nu]\,d\mathcal{H}^{N-1} .
	\]
	On the other hand \eqref{eqweakform3} gives $u\,[\z,\nu]=-|u|$ $\mathcal{H}^{N-1}$-a.e.
	on $\partial\Omega$, so that
	\begin{equation}\label{eq:energy-identity}
		\int_\Omega (\z,Du)+\int_{\partial\Omega}|u|\,d\mathcal{H}^{N-1}=\int_\Omega f u\,dx .
	\end{equation}
	We now estimate the left-hand side of \eqref{eq:energy-identity} from below. To this aim observe that it is not difficult to show that it follows from \eqref{eqweakform2} (see also \eqref{YoungFenchelIdweak}-\eqref{YoungFenchelIdweak2}) that it holds
	$$(\z,Du) \ge \ \Phi(|D^au|)\,\mathcal{L}^N+|D^su| \ge \alpha\,|Du|,$$
	as measures in $\Omega$.
	Therefore, for the left-hand side of \eqref{eq:energy-identity} (recall $\alpha \le 1$), we obtain
	\begin{equation}\label{eq:lhs}
		\int_\Omega (\z,Du)+\int_{\partial\Omega}|u|\,d\mathcal{H}^{N-1}
		\ \ge\ \alpha\left(\int_\Omega|Du|+\int_{\partial\Omega}|u|\,d\mathcal{H}^{N-1}\right)
		=\alpha\,\|u\|_{BV(\Omega)} .
	\end{equation}
	As for the right-hand side of \eqref{eq:energy-identity}, we apply the H\"older and the Sobolev inequality in order to deduce
	\begin{equation}\label{eq:rhs}
		\int_\Omega f u\,dx\le\|f\|_{L^{N,\infty}(\Omega)}\|u\|_{L^{1^*,1}(\Omega)}
		\le \tilde{S}_1\|f\|_{L^{N,\infty}(\Omega)}\,\|u\|_{BV(\Omega)} .
	\end{equation}
	Gathering \eqref{eq:energy-identity}, \eqref{eq:lhs} and \eqref{eq:rhs} we deduce
	\[
	\Big(\alpha-\tilde{S}_1\|f\|_{L^{N,\infty}(\Omega)}\Big)\,\|u\|_{BV(\Omega)}\le 0,
	\]
	which implies that $u\equiv0$ in $\Omega$.
\end{proof}

\begin{remark}
\rm We remark that the conclusion of Theorem \ref{thm:all-solutions-trivial} still holds once $\|f\|_{L^N(\Omega)}<\alpha\mathcal{S}_1^{-1}$,  i.e.  any solution is necessarily trivial. The proof is analogous to the one presented for the above theorem, employing the standard Sobolev inequality \eqref{sob} in place of the inequality \eqref{soblorentz}.

Finally observe that the   result of Theorem \ref{thm:all-solutions-trivial} is consistent with the known results for $\alpha=0$ (e.g. minimal surface operator, where $u$ is non-trivial \cite{gop}) and $\alpha=1$ (e.g. $1$-Laplacian, vanishing solutions \cite{CT,MST}).
\end{remark}

\begin{remark}
	\rm
	Theorem \ref{thm:all-solutions-trivial} shows that any solution to \eqref{P1} is trivial if $\|f\|_{L^{N,\infty}(\Omega)}$ is smaller than $\alpha\tilde{S}_1^{-1}$, which leads one  to wonder if solutions remain trivial as long as the norm is less than $\tilde{S}_1^{-1}$, as it happens for the $1$-Laplacian. 
	Here we aim to provide an argument suggesting that this is not the case.
	To begin with observe that
	$$
	\|f\|_{W^{-1,\infty}(\Omega)} = \sup \left\{ \langle f, v \rangle : v \in W_0^{1,1}(\Omega), \ \|v\|_{W_0^{1,1}(\Omega)} \le 1 \right\},
	$$
	then, reasoning as in the proof of Theorem \ref{thm:all-solutions-trivial},  one gets that the solutions are trivial once 	$\|f\|_{W^{-1,\infty}(\Omega)}<\alpha$.
	On the other hand, assume now that $\|f\|_{W^{-1,\infty}(\Omega)}>\alpha$. First, it is easy to check 
	$$
	\|f\|_{W^{-1,\infty}(\Omega)} \le \|\z\|_{\infty},
	$$
	holds true (in the case of the $1$-Laplacian (see \cite{MST}), we obtain the identity 
	$	\|f\|_{W^{-1,\infty}(\Omega)} = \|\z\|_{\infty}$).
	Notice that a sufficient condition to have $u \not\equiv 0$ is that $(\z, Du)$ is a positive measure, which is guaranteed if $\|\z\|_{\infty} > \alpha$.
	Indeed one has that $(u,\z)$ satisfies
	$$
	(\z, Du) = \Phi(|Du|) + \Phi^*(|\z|)\mathcal{L}^N \ge \Phi^*(|\z|)\mathcal{L}^N
	$$
	and $\Phi^*(t) > 0$ if and only if $|t| > \alpha$. This shows that non-trivial solutions appear in case
	$$
      \|\z\|_{\infty}\ge \|f\|_{W^{-1,\infty}(\Omega)}  > \alpha.
	$$
	It is also worth mentioning that it trivially holds that 
	$$
	\|f\|_{W^{-1,\infty}(\Omega)} \le  \tilde{S}_1\|f\|_{L^{N,\infty}(\Omega)},
	$$ 
	which completes our argument.
\end{remark}

\noindent{\bf Acknowledgments: }  The authors sincerely thank the anonymous referees for their careful reading and constructive suggestions, which have significantly improved the manuscript.
\\
F. Oliva and F. Petitta have been  partially supported by the Gruppo Nazionale per l’Analisi Matematica, la Probabilità e le loro Applicazioni (GNAMPA) of the Istituto Nazionale di Alta Matematica (INdAM).
\\
Marcos T.O. Pimenta is partially supported by FAPESP 2023/05300-4, 2022/16407-1 and CNPq 303868/2024-4, Brazil.
\\
The research of S. Segura de Le\'on is partially supported by Grant PID2022-136589NB-I00 funded
by MCIN/AEI/10.13039/501100011033 as well as by Grant RED2022-134784-T funded
by MCIN/AEI/10.13039/501100011033.


\begin{thebibliography}{99}

\bibitem{Adams}{ {\sc R.A. Adams}, {\em Sobolev spaces}, 2nd edition, Academic Press, New York (2003).}

\bibitem{AlvesSantosZhou} {{\sc J. Abrantes Santos, C.O. Alves and J. Zhou}, {\em Positive solutions for a class semipositone quasilinear problem with Orlicz-Sobolev critical growth}, Math. Nachr., 296, 4686--4711 (2023).}

\bibitem{SantosSoares} {\sc J. Abrantes Santos, S.H.M. Soares}, {\em Optimal design problems for a degenerate operator in Orlicz-Sobolev spaces}, Calc. Var. Partial Differential Equations, 428, 1035--1053 (2015).

\bibitem{AmbrosioFuscoPallara}{ {\sc L. Ambrosio, N. Fusco and D. Pallara}, {\em Functions of bounded variation and free discontinuity problems}, Oxford University Press., Oxford (2000).}

\bibitem{AndreuBallesterCasellesMazon} {{\sc F. Andreu, C. Ballester, V. Caselles and J. M. Maz\'on}, {\em Minimizing total variation flow}, C. R. Acad. Sci., Paris, Sr. I, Math., 331, No. 11, 867--872 (2000).}

\bibitem{AndreuBallesterCasellesMazonJFA} {{\sc F. Andreu, C. Ballester, V. Caselles and J. M. Maz\'on}, {\em The Dirichlet problem for the total variation flow}, J. Funct. Anal., 180, No. 2, 347--403 (2001).}

\bibitem{AndreuCasellesMazon} {}{{\sc F. Andreu, V. Caselles and J. M. Maz\'on}, {\em 
A parabolic quasilinear problem for linear growth functionals}, Rev. Mat. Iberoam. 18, No. 1, 135--185 (2002).}

\bibitem{An} {{\sc G. Anzellotti}, {\em Pairings between measures and bounded functions and compensated compactness}, Ann. Mat. Pura Appl. 135 (4),  293--318 (1983).}

\bibitem{anzellotti} {{\sc G. Anzellotti}, {\em The Euler equation for functionals with linear growth}, Trans. Amer. Math. Soc., 290, No. 2, 483--501 (1985).}

	\bibitem{BS} {{\sc C. Bennett and R. Sharpley}, {\em Interpolations of Operators}, Academic Press (1988).}
	
	\bibitem{Br} {{\sc H. Brezis}, {\em Functional analysis, Sobolev spaces and partial differential equations}, Springer (2011).}

\bibitem{CarvalhoGoncalvesSilva} {{\sc M. L. M. Carvalho, J. A. Gon\c{c}alves and E. D. Silva},  {\em On quasilinear elliptic problems without the Ambrosetti--Rabinowitz condition}, J. Math. anal. Appl., 426, 466--483 (2015).}


\bibitem{chl} {\sc I. Chlebicka},  {\em A pocket guide to nonlinear differential equations in Musielak-Orlicz spaces}, Nonlinear Anal. 175 1--27 (2018).

\bibitem{ClementHuidobroManasevichSchmitt} {{\sc P. Cl\'ement, M. Garc\'ia-Huidobro, R. Man\'asevich, K. Schmitt}, {\em Mountain pass type solutions for quasilinear elliptic equations}, Calc. Var. Partial Differential Equations, 11, 33--62 (2000).}


\bibitem{CT} {\sc M. Cicalese, C. Trombetti}, Asymptotic behaviour of solutions to $p$‐Laplacian equation, Asymp. Anal. 35 (1), 27-40, (2003)

	\bibitem{CF} {{\sc G. Q. Chen and H. Frid}, {\em Divergence-measure fields and hyperbolic conservation laws}, Arch. Ration. Mech. Anal. 147 (2),  89--118 (1999).}
	
	\bibitem{D} {{\sc F. Demengel}, {\em On some nonlinear partial differential equations involving the ``1''--Laplacian and critical Sobolev exponent}, ESAIM, Control Optim. Calc. Var. 4, 667--686 (1999).}
	
\bibitem{dt} {{\sc F. Demengel and R.  Temam}, {\em Convex functions of a measure and applications}, Indiana Univ. Math. J. 33 (5), 673--709 (1984).}

\bibitem{Donaldson} {{\sc T.K. Donaldson}, {\em Nonlinear elliptic boundary value problems in Orlicz–Sobolev spaces}, J. Differential Equations, 10, 507--528 (1971).}

\bibitem{EkelandTeman}{ {\sc I. Ekeland and R. Temam}, {\em Convex analysis and variational problems}, SIAM's Classics in Applied Mathematics, SIAM, New York (1999).}


\bibitem{HuidobroLeManasevichSchmitt} {{\sc M. Garc\'ia-Huidobro, V.K. Le, R. Man\'asevich, K. Schmitt}, {\em On principal eigenvalues for quasilinear elliptic differential operators: an Orlicz–Sobolev space setting}, NoDEA Nonlinear Differential Equations Appl., 6, 207--225 (1999).}

\bibitem{gop} {{\sc D. Giachetti, F. Oliva, F. Petitta}, {\em Bounded solutions for non-parametric mean curvature problems with nonlinear terms}, J. Geom. Anal. 34, No. 9, Paper No. 265, 32 p. (2024). }

\bibitem{gia} {\sc M. Giaquinta}, {\em On the Dirichlet problem for surfaces of prescribed mean curvature}, Manuscripta Math. 12, 73--86 (1974).

\bibitem{g76} {\sc E. Giusti}, {\em Boundary value problems for non-parametric  surfaces of prescribed mean curvature}, Ann. Scuola Norm. Sup. Pisa Cl. Sci. (4) 3 (3), 501--548 (1976).

\bibitem{g78} {\sc E. Giusti},   {\em On the Equation of Surfaces of Prescribed Mean Curvature. Existence and uniqueness without
boundary conditions}, Invent. math. 46, 111--137 (1978).

\bibitem{Gossez} {{\sc J.P. Gossez}, {\em Nonlinear elliptic boundary value problems for equations with rapidly (or slowly) increasing coefficients}, Trans. Amer. Math. Soc., 190, 163--205 (1974).}

\bibitem{GossezManasevich} {{\sc J.P. Gossez and R. Man\'asevich}, {\em On a nonlinear eigenvalue problem in Orlicz–Sobolev spaces}, Proc. Roy. Soc. Edinburgh Sect. A, 132, 891--909 (2002).}


\bibitem{K1} {\sc B. Kawohl}, { On a family of torsional creep problems}, J. Reine Angew. Math. 410, 1-22, (1990).

\bibitem{K2} {\sc B. Kawohl}, {From $p$-Laplace to mean curvature operator and related questions}, In: Progress in Partial Differential Equations: the Metz Surveys, Pitman Res. Notes Math. Ser. 249, 40-56, (1991).

\bibitem{LL} {{\sc J. Leray and J.L. Lions}, {\em Quelques reesultats de Vis\'ik sur le problemes elliptiques non lineaires par le methodes de Minty-Browder}, Bul. Soc. Math. France 93, 97--107 (1965).}

%\bibitem{MazonNA2016} {{\sc J. M. Maz\'on}, {\em The Euler-Lagrange equation for the anisotropic least gradient problem}, Nonlinear Anal. Real World Appl., 31, 452--472 (2016).}

\bibitem{MRST} {\sc A. Mercaldo, J.D. Rossi, S. Segura de Le\'on and C. Trombetti}, {\em Behaviour of $p$-Laplacian problems with Neumann boundary conditions when $p$ goes to 1},
Commun. Pure Appl. Anal. 12, No. 1, 253-267 (2013). 

\bibitem{MST} {\sc A. Mercaldo, S. Segura de Le\'on and C. Trombetti}, {\em On the behaviour of the solutions to p--Laplacian equations as p goes to 1}, Publ. Mat., Barc. 52, No. 2, 377-411 (2008).

\bibitem{MihailescuRepovs} {\sc M. Mih\u{a}ilescu and D. Repov\v{s}}, {\em Multiple solutions for a nonlinear and non-homogeneous problem in Orlicz-Sobolev spaces}, Appl. Math. Comput., 217, 6624--6632 (2011).

\bibitem{MontenegroPimenta} {\sc M. Montenegro and M.T.O. Pimenta}, {\em Existence of a BV solution for a mean curvature equation}, J. Differential Equations, 299, 51--64 (2021).

\bibitem{OlivaPetittaSegura} {\sc F. Oliva, F. Petitta and S. Segura de Le\'on}, {\em The role of absorption terms in Dirichlet problems for the prescribed mean curvature equation}, NoDEA Nonlinear Differential Equations Appl., 31, 53 (2024).

	\bibitem{PKF}  {\sc L. Pick, A. Kufner, J. Oldrich and S. Fuc\'ik},  Function Spaces, 1. Berlin, Boston: De Gruyter, 2012.

\bibitem{stevsch}{\sc  C. Scheven, T. Schmidt}, {\em  BV supersolutions to equations of 1-Laplace and minimal surface type}, J. Differ. Equations 261, (2016) 1904--1932.

\bibitem{SchevenSchmidtAIHP} {{\sc C. Scheven and T. Schmidt},  {\em On the dual formulation of obstacle problems for the total variation and the area functional}. Ann. I. H. Poincar\'e Anal. Non Lin\'eaire 35, 1175--1207 (2018).}

%\bibitem{silva} {\sc E.D. Silva, M.L.  Carvalho, and  J.C. de Albuquerque},  {\em Revised regularity results for quasilinear elliptic problems driven by the 
%-Laplacian operator}, manuscripta math. 161, 563-582 (2020).

\bibitem{Stampacchia} {{\sc G. Stampacchia}, {\em Le probl\`eme de Dirichlet pour les \'equations elliptiques du seconde ordre  \`a coefficientes discontinus},  Ann.
Inst. Fourier (Grenoble) 15, 189--258  (1965).}
 
\bibitem{Te} {{\sc R. Temam}, {\em Solutions g\'en\'eralis\'ees de certaines \'equations du type hypersurfaces minima},  Arch. Ration. Mech. Anal. 44, 121--156 (1971).}

%\bibitem{r} {{\sc Yu. G. Reschetniak}, {\em Weak convergence of completely additive vector functions on a set, Sibirsk. Mat. Z.}, 9,1386--1394 (1968) (Translated)}

\bibitem{Z} {\sc A.C. Zaanen}, {\em  Linear analysis} 3rd edition. 
Bibliotheca Mathematica. Vol. 2. Groningen: P. Noordhoff Ltd. (1960). 

	\bibitem{Zh} {{\sc X. Zhou}, {\em An evolution problem for plastic antiplanar shear},  Appl. Math. Optimization 25, No. 3, 263--285 (1992).}

\end{thebibliography}
\end{document}